\documentclass[MSNbibl,number,citesort,seceqn,dvips]{arxbj}

\aid{0}
\volume{19}
\issue{5B}
\pubyear{2013}
\firstpage{2294}
\lastpage{2329}
\doi{10.3150/12-BEJ453} 

\makeatletter
\newcommand{\rrVert}{\Vert}
\newcommand{\rrvert}{\vert}
\newcommand{\llVert}{\Vert}
\newcommand{\llvert}{\vert}
\newtheorem{theo}{Theorem}[section]
\newtheorem{lem}{Lemma}[section]
\newtheorem{prop}{Proposition}[section]
\newtheorem{cor}{Corollary}[section]
\newproclaim{dfn}{Definition}[section]
\newremark{ex}{Example}[section]
\newremark{rem}{Remark}[section]
\newcommand{\Cov}{\operatorname{Cov}}
\newcommand{\Lip}{\operatorname{Lip}}
\newcommand{\Corr}{\operatorname{Corr}}
\newcommand{\Dom}{\operatorname{Dom}}
\def\XX{\mathbb{X}}
\def\Bf{\mathcal{B}}
\def\Pf{\mathcal{P}}
\def\Af{\mathcal{A}}
\def\Ff{\mathcal{F}}
\renewcommand{\Re}{\mathbb{R}}
\def\eps{\varepsilon}
\def\1{\mathbf{1}}

\newcommand{\eqref}[1]{(\ref{#1})}
\makeatother

\begin{document}
\begin{frontmatter}

\title{Ergodicity and mixing bounds for the Fisher--Snedecor diffusion}
\runtitle{Ergodicity and mixing bounds for the Fisher--Snedecor diffusion}

\begin{aug}
\author[1]{\fnms{A.M.} \snm{Kulik}\corref{}\thanksref{1}\ead[label=e1]{kulik@imath.kiev.ua}} \and
\author[2]{\fnms{N.N.} \snm{Leonenko}\thanksref{2}\ead[label=e2]{LeonenkoN@Cardiff.ac.uk}}
\runauthor{A.M. Kulik and N.N. Leonenko} 
\address[1]{Institute of Mathematics, Ukrainian National Academy of
Sciences, 01601 Tereshchenkivska str. 3, Kyiv, Ukraine. \printead{e1}}
\address[2]{School of Mathematics, Cardiff University, Senghennydd
Road, Cardiff CF24 4AG, UK.\\ \printead{e2}}
\end{aug}

\received{\smonth{10} \syear{2011}}
\revised{\smonth{5} \syear{2012}}

%
\begin{abstract}
We consider the Fisher--Snedecor diffusion; that is, the
Kolmogorov--Pearson diffusion with the Fisher--Snedecor invariant
distribution. In the nonstationary setting, we give explicit
quantitative rates for the convergence rate of respective
finite-dimensional distributions to that of the
stationary Fisher--Snedecor diffusion, and for {the $\beta$-mixing
coefficient} of this diffusion. As an application, we prove the law of
large numbers and the central limit theorem for additive functionals of
the Fisher--Snedecor diffusion and construct $P$-consistent and asymptotically
normal estimators for the parameters of this diffusion given its
nonstationary observation.
\end{abstract}

%

\end{frontmatter}

\section{Introduction}

In this paper, we investigate the Markov process $X$, valued in
$(0,\infty)$%
, defined by the nonlinear stochastic differential equation
%
\begin{equation}
\label{fssde1} \mathrm{d}X_{t} = -\theta (X_{t} - \kappa )
\,\mathrm{d}t + \sqrt{2 \theta X_{t} \biggl( \frac{ X_{t}}{\beta/2 - 1} +
\frac{ \kappa}{\alpha/2 } \biggr)} \,\mathrm{d}W_{t},\qquad t \geq0.
\end{equation}
Such a process belongs to the class of diffusion processes with invariant
distributions from the \emph{Pearson family}, introduced by K. Pearson
\cite{Pearson} in 1914 in order to unify some of the most important statistical
distributions. The study of such processes was started in the 1930s by
A.N. Kolmogorov \cite{Kolmogorov,Shiryayev}, hence it seems
appropriate to call this important class of processes the \emph{%
Kolmogorov--Pearson (KP) diffusions}. For a more detailed discussion of KP
diffusions, we refer to recent papers \cite{FormanSorensen,ShawMunir} and \cite{ALS1}.

When $\alpha, \beta> 2$, the process $X$ defined by (\ref{fssde1}) is
ergodic \cite{Genon-Catalot}. Under the particular choice $\kappa={%
\beta/(\beta-2)}$, respective unique invariant distribution coincides with
the \emph{Fisher--Snedecor distribution} $\mathcal{F S}(\alpha, \beta)$
with $%
\alpha,\beta$ degrees of freedom; that is, its probability density is given
by
%
\begin{equation}
\label{fsdens} \mathrm{\mathfrak{fs}}(x) ={\frac{1}{x B(\alpha/2, \beta
/2)}} \biggl(%
\frac{\alpha x}{\alpha x+\beta} \biggr)^{\alpha/2} \biggl(\frac
{\beta
}{\alpha
x+ \beta}
\biggr)^{\beta/2},\qquad x>0.
\end{equation}
This is the reason to call the process $X$ defined by (\ref{fssde1})
\emph{%
the Fisher--Snedecor diffusion}. Together with \emph{the reciprocal gamma}
and \emph{the Student} diffusions, the Fisher--Snedecor diffusion forms the
class of the so-called \emph{heavy-tailed} KP diffusions. Statistical
inference for three heavy-tailed KP diffusions is developed in the recent
papers \cite{Leonenko1,Leonenko2} and \cite{ALS1} in the situation
where the stationary version of the respective diffusion is observed.

In this paper, we consider the Fisher--Snedecor diffusion (\ref{fssde1}) in
the nonstationary setting; that is, with arbitrary distribution of the
initial value $X_0$. We give explicit quantative rates for \emph{the
convergence rate} of respective finite-dimensional distributions to
that of
the stationary Fisher--Snedecor diffusion, and for \emph{the $\beta$-mixing
coefficient} of this diffusion. Same problems for the reciprocal gamma and
the Student diffusions were considered in \cite{AV1} and \cite{AV2},
respectively. Similarly to \cite{AV1} and \cite{AV2}, our way to
treat this
problem is based on the general theory developed for (possibly nonsymmetric
and nonstationary) Markov processes, although there is a substantial
novelty in the form taken by the \emph{Lyapunov-type condition}
(typical in
the field) in our setting.

As an application, we prove the law of large numbers (LLN) and the central
limit theorem (CLT) for additive functionals of the Fisher--Snedecor
diffusion. Note that, for the stationary version of the diffusion, these
limit theorems are well known: LLN is provided by the Birkhoff--Khinchin
theorem, and CLT is available either in the form based on the $\alpha$%
-mixing coefficient of a stationary sequence or process (see \cite{Hall}),
or in the form formulated in terms of the $L_2$-semigroup associated with
the Markov process (see \cite{Bhatt}). Our considerations are based on the
natural idea to extend these results to the nonstationary setting
using the
bounds for the deviation between the stationary and nonstationary versions
of the process. The way we carry out this idea differs, for instance, from
those proposed in \cite{Bhatt}, Theorem 2.6, or in \cite{Shiryayev2},
Section %
4.II.1.10, and is based on the notion of an \emph{(exponential) $\phi
$-coupling%
}, introduced in \cite{Kulik} as a tool for studying convergence rates
of $%
L_p$-semigroups, generated by a Markov process, and spectral properties of
respective generators.

The modified version of the Lyapunov-type condition, mentioned above,
implies a substantial difference between the asymptotic properties of the
finite-dimensional distributions themselves and their continuous-time
averages, see Theorem \ref{t37} and Remark \ref{r32} below. An important
consequence is that, in the continuous-time version of our CLT, the
observable functional \emph{may fail to be square integrable} w.r.t. the
invariant distribution of the process. This interesting effect
seemingly has
not been observed in the literature before.

Finally, we apply the above results and provide a statistical analysis for
the Fisher--Snedecor diffusion. In the situation where a nonstationary
version of the diffusion $X$ is observed, we prove that respective empirical
moments and empirical covariances are $P$-consistent, asymptotically normal,
and (under some additional assumptions on the initial distribution of $X$)
asymptotically unbiased. Then, using the method of moments, we
construct $P$%
-consistent and asymptotically normal estimators for the parameter
$(\alpha,
\beta, \kappa, \theta)$ given either the discrete-time or the
continuous-time observations of a nonstationary version of the
Fisher--Snedecor diffusion. To keep the current paper reasonably short, we
postpone the explicit calculation of the asymptotic covariance matrices and
a more detailed discussion of other statistical aspects to the subsequent
paper \cite{KulikLeonenko}.

\section{Preliminaries}

In this section, we introduce briefly main objects, assumptions, and
notation.

For the Fisher--Snedecor diffusion \eqref{fssde1}, the drift
coefficient $a(x)
$ and the diffusion coefficient $\sigma(x)$ are respectively given by
%
\begin{equation}
\label{dd} a(x) = -\theta (x - \kappa ),\qquad \sigma(x) = \sqrt{2 \theta x
\biggl( \frac{x}{\beta/2 - 1} +\frac{ \kappa}{\alpha/2 } \biggr)},
\end{equation}
and our standing assumptions on the parameters are
%
\begin{equation}
\label{parameters} \theta>0,\qquad \kappa>0,\qquad \beta>2,\qquad \alpha> 2.
\end{equation}
We assume that, on a proper probability space $(\Omega, P, \mathcal
{F})$, independent Wiener process $W$ and random variable $X_0$ taking
values in $(0, \infty)$ are well defined. Then, because the
coefficients (\ref{dd}) are continuously differentiable inside $
(0,\infty)$, the unique
strong solution to equation \eqref{fssde1} with the initial condition
$X_0$ is well defined up to the random
time moment $T_{0,\infty}$ of its exit from $(0,\infty)$.

For $x\in(0,\infty)$, the
corresponding \emph{scale density} equals
%
\begin{equation}
\mathfrak{s}(x) =\exp \biggl( -\int_{1}^x
\frac{2 a(u)}{\sigma^{2}(u)}\, \mathrm{d}u \biggr)= C x^{-\alpha/2} \biggl(x + {
\frac{\kappa(\beta
-2)}{\alpha}}%
\biggr)^{\alpha/2+\beta/2 -1}. \label{scale}
\end{equation}
Here and below, by $C$ we denote a constant, which can be (but is not)
expressed explicitly; the value of $C$ can vary from place to place. It
follows from the standing assumption (\ref{parameters}) that
\[
\int_{x}^{\infty} \mathfrak{s}(y)\,\mathrm{d}y =
\infty,\qquad \int_{0}^{x} \mathfrak{s}(y)\,
\mathrm{d}y = \infty,\qquad x\in(0, \infty),
\]
and consequently both $0$ and $\infty$ are
unattainable points for the diffusion $X$, that is, the random time moment
$T_{0,\infty}$ is a.s. infinite for any positive initial condition
$X_0$ (e.g., \cite{Karlin2}, Chapter 18.6). This means that (\ref{fssde1}) uniquely determines a time-homogeneous strong Markov process $X$
with the state space $\mathbb{X}=(0,\infty)$. In the sequel, we
consider $%
\mathbb{X}$ as a locally compact metric space with the metric $%
d(x,y)=|x-y|+|x^{-1}-y^{-1}|.$

Let us introduce the notation. By $P_t(x,\mathrm{d}y)$, we denote the transition
probabilities of the process $X$. By $\mathcal{P}$ we denote the class of
probability distributions on the Borel $\sigma$-algebra on $\mathbb
{X}$. For
any $\mu\in\mathcal{P}$, we denote by $P_\mu$ the distribution in
$C(\Re^+,
\mathbb{X})$ of the solution to (\ref{fssde1}) with the distribution
of $X_0$
equal $\mu$, and write $E_\mu$ for the respective expectation. When $%
\mu=\delta_x$, the measure concentrated at the point $x\in\mathbb
{X}$, we
write $P_x,E_x$ instead of $P_\mu, E_\mu$. For any $\mu\in\mathcal
{P}$ we
denote by $\mu_{t_1,\ldots, t_m}, 0\leq t_1<\cdots<t_m, m\geq1$ the
family of
finite-dimensional distributions of the process $X$ with the initial
distribution $\mu$; that is,
%
\begin{eqnarray}
\label{mut} %
\mu_{t_1, \ldots,
t_m}(A)&=&\int_\XX
\int_{A}P_{t_1}(x,\mathrm{d}x_1)P_{t_2-t_1}(x_1,
\mathrm{d}x_2)\cdots P_{t_m-t_{m-1}}(x_{m-1},
\mathrm{d}x_m)\mu(\mathrm{d}x)
\nonumber
\\[-8pt]
\\[-8pt]
&=&P_\mu\bigl((X_{t_1}, \ldots, X_{t_m})\in A
\bigr),\qquad A\in\mathcal{B}\bigl({\Bbb X}^m\bigr).
\nonumber
\end{eqnarray}

By $\mathbb{F}^X=\{\mathcal{F}_t^X, t\geq0\}$, we denote the natural
filtration of the process $X$. A measurable function $f\dvtx \mathbb{X}\to
\Re$
is said to belong to the domain of the \emph{extended generator
$\mathcal{A}$%
} of the process $X$ if there exists a measurable function $g\dvtx \mathbb
{X}\to
\Re$ such that the process
\[
f(X_t)-\int_0^tg(X_s)
\,\mathrm{d}s,\qquad t\in\Re^+
\]
is well defined and is an $\mathbb{F}^X$-martingale w.r.t. to any measure
$P_x, x\in\mathbb{X}$. For such a pair $(f,g)$, we write $f\in
\Dom(\mathcal{A%
})$ and $\mathcal{A} f= g$.

For a measurable function $\phi\dvtx \mathbb{X}\to[1,\infty)$ and a signed
measure $\varkappa$ on $\mathcal{B}(\mathbb{X}^m)$, define the \emph
{%
weighted total variation norm}
\[
\|\varkappa\|_{\phi,\mathrm{var}}=\int_{\mathbb{X}^m} \bigl(
\phi(x_1)+\cdots +\phi (x_m)%
\bigr)|\varkappa|(
\mathrm{d}x),
\]
where $|\varkappa|=\varkappa^++\varkappa^-$ and $\varkappa
=\varkappa^+-%
\varkappa^-$ is the Hahn decomposition of $\varkappa$. Frequently, we will
use functions $\phi$ of the form
%
\begin{equation}
\label{phidef} \phi=\phi_\lozenge+\phi_\blacklozenge,
\end{equation}
where $\phi\geq1$, $\phi_\lozenge,\phi_\blacklozenge\in
C^2(0,\infty
)$, $\phi_\lozenge=0$
on $[2,\infty)$, $\phi_\blacklozenge=0$ on $(0,1]$,
\[
\phi_\lozenge(x)=x^{-\gamma}\qquad\mbox{for $x$ small enough},\qquad
\phi_\blacklozenge(x)=x^{\delta}\qquad\mbox{for $x$ large enough}
\]
with nonnegative $\gamma, \delta$.

The \emph{$\beta$-mixing} (or \emph{complete regularity}, or \emph{the
Kolmogorov}) coefficient is defined as
%
\begin{equation}
\label{beta-m} \beta^\mu(t)=\sup_{s\geq0}E_\mu
\sup_{B\in\mathcal{F}^X_{\geq
t+s}}\bigl|P_\mu\bigl(B|%
\mathcal{F}^X_s
\bigr)-P_\mu(B)\bigr|,\qquad\mu\in\mathcal{P}, t\in\Re^+,
\end{equation}
where $\mathcal{F}^X_{\geq r}$ for a given $r\geq0$ denotes the
$\sigma
$%
-algebra generated by the values of the process $X$ at the time moments
$%
v\geq r$. In particular, the state-dependent $\beta$-mixing
coefficient is
defined by
%
\begin{equation}
\label{beta-m1} \beta_x(t)=\sup_{s\geq0}E_x
\sup_{B\in\mathcal{F}^X_{\geq
t+s}}\bigl|P_x\bigl(B|%
\mathcal{F}^X_s
\bigr)-P_x(B)\bigr|,\qquad x\in\mathbb{X}, t\in\Re^+
\end{equation}
(in this case, the initial distribution $\mu=\delta_x$), and the
stationary $%
\beta$-mixing coefficient is defined by
%
\begin{equation}
\label{beta-m2} \beta(t)=\sup_{s\geq0}E_\pi
\sup_{B\in\mathcal{F}^X_{\geq
t+s}}\bigl|P_\pi \bigl(B|%
\mathcal{F}^X_s
\bigr)-P_\pi(B)\bigr|,\qquad x\in\mathbb{X}, t\in\Re^+;
\end{equation}
here and below, $\pi$ denotes the (unique) invariant distribution for the
process $X$. For more information about various types of mixing coefficients
see, for example, \cite{Bradley}.

\section{Main results}\label{s2}

Here, we formulate the main results of the paper. The proofs are postponed
to Section~\ref{s4}.

\subsection{Distributional properties of the Fisher--Snedecor diffusion}

The following two basic properties of the Fisher--Snedecor diffusion
will be
used in the further analysis of its ergodic behavior.

\begin{prop}\label{p21}
1. (Lyapunov-type condition). Let $\phi$ to have the form
(\ref{phidef}) with
%
\begin{equation}
\label{gammadelta} \gamma<{\frac{\alpha}{2}}-1,\qquad \delta<{
\frac{\beta}{2}}.
\end{equation}

Then $\phi\in \Dom(\mathcal{A})$ and
%
\begin{equation}
\label{Af} \mathcal{A} \phi=a\phi^{\prime}+{\frac{1}{2}}
\sigma^2\phi^{\prime\prime}.
\end{equation}
In addition, there exist a segment $[u,v]\subset(0,\infty)$ and positive
constants $c,C$ such that
%
\begin{equation}
\label{lyap} \mathcal{A} \phi(x)\leq-c\phi(x)+C\1_{[u,v]}(x).
\end{equation}

2. (Local minorization condition). For every segment $[u,v]\subset
\mathbb{X}$
there exist $T>0$, another segment $[u^{\prime},v^{\prime}]\subset
\mathbb{X}
$ and a constant $c_{u,v,u^{\prime},v^{\prime},T}>0$ such that for
every $%
x\in[u,v]$ and every Borel set $A\subset[u^{\prime},v^{\prime}]$
\[
P_T(x,A)\geq c_{u,v,u^{\prime},v^{\prime},T}\int_A
\mathrm{d}y.
\]
\end{prop}

The following moment bound is a well known corollary of the Lyapunov-type
condition (see, e.g., Section 3.2 in \cite{K09} and references therein).

\begin{cor}
\label{c31} In the conditions and notation of statement 1 in
Proposition \ref{p21}, we have
\[
\int_\XX\phi \,\mathrm{d}\mu_t\leq{
\frac{C}{c}}+\mathrm{e}^{-c t}\int_\XX\phi\,
\mathrm{d}\mu, \qquad t\in\Re^+.
\]
In addition, there exists an invariant measure $\mu^*\in\mathcal{P}$ such
that
\[
\int_\XX\phi \,\mathrm{d}\mu^*<+\infty.
\]
\end{cor}

Because the Fisher--Snedecor diffusion is ergodic, the latter statement can
be interpreted as the following fact about its (unique) invariant
distribution $\pi$:
%
\begin{equation}
\label{moments} \int_\XX x^{-\gamma}\pi(
\mathrm{d}x)<+\infty,\qquad \int_\XX x^{\delta}\pi (
\mathrm{d}x)<+\infty
\end{equation}
as soon as positive $\gamma, \delta$ satisfy (\ref{gammadelta}). On the
other hand, the probability density $\mathfrak{p}$ of the invariant
distribution $\pi$ is proportional to $\sigma^{-2}\mathfrak{s}^{-1}$ (e.g.,
see \cite{ALS1}), and straightforward calculation shows that (\ref{moments})
holds true if, and only if,
%
\begin{equation}
\label{gammadelta1} \gamma<{\frac{\alpha}{2}},\qquad \delta<{
\frac{\beta}{2}}.
\end{equation}
Clearly, the first bound in (\ref{gammadelta1}) is weaker than the one
in (%
\ref{gammadelta}). Such a discrepancy indicates that, in the current
setting, the Lyapunov-type condition (\ref{lyap}) is not precise, in a
sense. This observation motivates the following extension of the above
results. Define the family of \emph{Ces\`{a}ro means} of finite-dimensional
distributions of $X$ by
%
\begin{equation}
\label{cesaro} \mu^t_{t_1, \ldots, t_m}={\frac{1}{t}}\int
_0^t\mu_{t_1+s, \ldots,
t_m+s} \,\mathrm{d}s,\qquad
t>0,0\leq t_1<\cdots<t_m, m\geq1.
\end{equation}

\begin{prop}
\label{p32} 1. (Modified Lyapunov-type condition). Let $\phi$ have the
form (%
\ref{phidef}) with positive $\gamma, \delta$ satisfying (\ref
{gammadelta1}%
). Then there exists a nonnegative function $\psi\in \Dom(\mathcal{A})$,
satisfying (\ref{lyap}) and such that
%
\begin{equation}
\label{est} \mathcal{A} \psi\leq-c^{\prime}\phi^{1+\varepsilon}+C^{\prime}
\end{equation}
with some positive constants $c^{\prime},C^{\prime}, \varepsilon$.

2. (Moment bounds for Ces\`{a}ro means). In the conditions and notation of
statement 1, let $c, C$ be the constants from the relation (\ref
{lyap}) for
the function $\psi$. Then, for arbitrary $m\geq1, 0\leq t_1<\cdots< t_m$,
%
\begin{equation}
\label{Cesarobounds} \int_{\XX^m} \bigl(\phi(x_1)+
\cdots+\phi(x_m) \bigr)^{1+\varepsilon
}\mu^t_{t_1,
\ldots, t_m}(
\mathrm{d}x)\leq m^\eps \biggl({\frac{C^{\prime}}{c^{\prime
}}}+{
\frac
{C}{%
cc^{\prime}t}}+\frac{1}{ c't} \biggr)\int_\XX\psi
\,\mathrm{d}\mu.
\end{equation}
\end{prop}

\begin{rem}
\label{r31} Let $\mu=\delta_x$, then (\ref{Cesarobounds}) with $m=1$
and $%
t_1=0$ yields
\[
\sup_{t\geq1}{\frac{1}{t}}\int_0^t
\int_\XX\phi \,\mathrm{d}\mu_s \,\mathrm{d}s<
\infty.
\]
On the other hand, by Theorem \ref{t31} below we have
\[
{\frac{1}{t}}\int_0^t
\mu_s \,\mathrm{d}s\Rightarrow\pi,\qquad t\to\infty.
\]
These two observations, combined with the proper version of the Fatoux lemma
(e.g., \cite{Billingsley}, Theorem 5.3) provide that $\phi$ is integrable
w.r.t $\pi$. This means that the moment bound (\ref{Cesarobounds})
yields (%
\ref{moments}) under (\ref{gammadelta1}), and hence resolves the discrepancy
discussed above.
\end{rem}

\subsection{\texorpdfstring{Coupling, ergodicity, and $\beta$-mixing}
{Coupling, ergodicity, and beta-mixing}}

This section collects the results about the ergodic behavior of the
Fisher--Snedecor diffusion. For our further needs, it will be convenient to
introduce explicitly and discuss separately the notion of an \emph{%
exponential $\phi$-coupling.}

By the common terminology, a \textit{coupling} for a pair of processes $U,V$
is any two-component process $Z=(Z^1,Z^2)$ such that $Z^1$ has the same
distribution with $U$ and $Z^2$ has the same distribution with $V$.
Following this terminology, for a Markov process $X$ and every $\mu
,\nu
\in%
\mathcal{P}$, we consider two versions $X^\mu,X^\nu$ of the process
$X$ with
the initial distributions equal {to} $\mu$ and $\nu$, respectively,
and call
\textit{$(\mu,\nu)$-coupling for the process $X$} any two-component
process $%
Z=(Z^1,Z^2)$ which is a coupling for $X^\mu, X^\nu$.

\begin{dfn}
\label{d13} The Markov process $X$ \emph{admits an exponential $\phi$-coupling} if there exists an invariant measure $\pi$ for this process and
positive constants $C,c$ such that, for every $\mu\in\mathcal{P}$, there
exists a $(\mu,\pi)$-coupling $Z=(Z^1,Z^2)$ with
%
\begin{equation}
\label{phi-coup} E \bigl[\phi\bigl(Z^1_t\bigr)+\phi
\bigl(Z_t^2\bigr) \bigr]\1_{Z_t^1\not
=Z_t^2}\leq C
\mathrm{e}^{-c t}\int_\XX\phi \,\mathrm{d}\mu,\qquad
t\geq0.
\end{equation}
\end{dfn}

The coupling construction is a traditional tool for proving the ergodicity.
In \cite{Kulik}, it was proposed to introduce a separate notion of an
exponential $\phi$-coupling, and it was demonstrated that such a notion
is a
convenient tool for studying convergence rates of $L_p$-semigroups,
generated by a Markov process, and spectral properties of respective
generators. In Section \ref{s44} below, we will see that this notion
is also
efficient for proving LLN and CLT. With this application in mind, we have
changed slightly Definition \ref{d13}, if to compare it with the one given
in \cite{Kulik}: here, we consider all probability measures $\mu\in
\mathcal{%
P}$ as possible initial distributions, while in \cite{Kulik} only measures
of the form $\mu=\delta_x, x\in\mathbb{X}$ are considered.

\begin{theo}
\label{t31} Let $\phi$ be defined by (\ref{phidef}) with $\gamma
,\delta$
satisfying (\ref{gammadelta}). Then the following statements hold.

1. The Fisher--Snedecor diffusion admits an exponential $\phi$-coupling.

2. Finite-dimensional distributions of the Fisher--Snedecor diffusion admit
the following convergence rate in the weighted total variation norm
with the
weight $\phi$: for any $m\geq1, 0\leq t_1<\cdots<t_m,$
%
\begin{equation}
\label{41} \|\mu_{t+t_1, \ldots, t+t_m}-\pi_{t_1, \ldots, t_m}\|_{\phi,\mathrm{var}}\leq m C
\mathrm{e}^{-c
t}\int_{\mathbb{X}}\phi \,\mathrm{d}\mu,\qquad
\mu\in\mathcal{P}, t\geq0.
\end{equation}
Here the constants $C,c$ are the same as in the bound (\ref{phi-coup}) in
the definition of an exponential $\phi$-coupling.

3. The Fisher--Snedecor diffusion admits the following bound for the
$\beta$%
-mixing coefficient:
%
\begin{equation}
\label{mixing} \beta^\mu(t)\leq C^{\prime}
\mathrm{e}^{-ct}\int_\XX\phi \,\mathrm{d}\mu,\qquad
\mu\in \mathcal{%
P}, t\geq0.
\end{equation}
Here the constant $c$ is the same as in the bound (\ref{phi-coup}), and
$%
C^{\prime}$ a positive constant, which can be given explicitly (see
(\ref{Cpr}) below).
\end{theo}

From (\ref{mixing}) and Corollary \ref{c31}, we get the following
bounds for
state-dependent and stationary $\beta$-mixing coefficients:
\begin{eqnarray*}
\beta_x(t)&\leq& C^{\prime}\mathrm{e}^{-ct}\phi(x),
\qquad x\in\mathbb{X}, t\geq 0,
\\
\beta(t)&\leq& C^{\prime\prime}\mathrm{e}^{-ct},\qquad t\geq0,
C^{\prime\prime}:=C^{\prime}\int_\XX\phi \,\mathrm{d}
\pi<+\infty.
\end{eqnarray*}

Note that the general theory for (possibly nonsymmetric and nonstationary)
Markov processes provides convergence rates like (\ref{41}), for
example, \cite%
{DMT95}, and bounds for $\beta$-mixing coefficients like (\ref{mixing}),
for example, \cite{Ver87}, under a proper combination of ``recurrence''
and ``local
irreducibility'' conditions. In our context, these conditions are
provided by
Proposition~\ref{p21}.

Apart with the convergence rate (\ref{41}), we give the following more
specific bound for continuous-time averages of the family $\{\mu_{t_1,
\ldots, t_m}\}$.

\begin{theo}
\label{t37} Let $\phi$ be defined by (\ref{phidef}) with $\gamma
,\delta$
satisfying (\ref{gammadelta1}), and $\psi$ be the function from
Proposition %
\ref{p32}.

Then for every $m\geq1$ there exists a constant $C_m$ such that
%
\begin{equation}
\label{45} \biggl\llVert \int_{0}^T (
\mu_{t+t_1, \ldots, t+t_m}-\pi_{t_1, \ldots, t_m}) \,\mathrm{d}t\biggr\rrVert_{\phi, \mathrm{var}}
\leq C_m\int_\XX\psi \,\mathrm{d}\mu,\qquad \mu
\in \mathcal{P}%
, T\geq0.
\end{equation}
\end{theo}

\begin{rem}
\label{r32} Clearly, (\ref{41}) provides a bound, similar to (\ref{45}),
with $\phi$ instead of $\psi$ in the right-hand side. This bound is weaker
than (\ref{45}) because $\psi(x)=\mathrm{o}(\phi(x))$ as $x\to0$ or $x\to
\infty$.
In addition, Theorem \ref{t37} requires (\ref{gammadelta1}), which is weaker
than respective assumption (\ref{gammadelta}) in Theorem \ref{t31}.
In this
sense, for continuous-time averages of the family $\{\mu_{t_1, \ldots,
t_m}\}$
Theorem \ref{t37} provides a substantially more precise information than
Theorem \ref{t31} does.
\end{rem}

\subsection{The law of large numbers and the central limit theorem}

In this section, we formulate LLN and CLT for additive functionals of the
Fisher--Snedecor diffusion $X$. Below, $X^{st}_t, t\in(-\infty, \infty)$
denotes the stationary version of $X$; that is, the strictly stationary
process such that for every $m\geq1$ and $t_1<\cdots<t_m$ the joint
distribution of $X^{st}_{t_1}, \ldots, X^{st}_{t_m}$ equals $\pi_{0,
t_2-t_1,\ldots, t_m-t_1}$ (heuristically, $X^{st}$ is ``a solution to
(\ref{fssde1}), which is defined on the whole time axis and starts at
$-\infty$ from the invariant distribution~$\pi$'').

We consider separately the
discrete-time and the continuous-time cases.

\begin{theo}[(Discrete-time case)]\label{t32}  Let, for some $r,k\geq1$, a vector-valued
function
\[
f=(f_1, \ldots, f_k)\dvtx \mathbb{X}^r\to
\Re^k
\]
be such that for any $i=1, \ldots, k$ for some $\gamma_i, \delta_i$
satisfying (\ref{gammadelta})
%
\begin{equation}
\label{fi} \bigl|f_i(x)\bigr|\leq C\sum_{j=1}^r
\bigl(x^{-\gamma_i}_j+x^{\delta_i}_j \bigr),
\qquad x=(x_1, \ldots, x_r)
\end{equation}
with some constant $C$.

Then the following statements hold true.

1. (LLN). For arbitrary initial distribution $\mu$ of $X$ and arbitrary
$%
t_1, \ldots, t_r\geq0$,
%
\begin{equation}
\label{LLNd} {\frac{1}{n}}\sum_{l=1}^nf
(X_{t_1+l}, \ldots, X_{t_r+l} )\to a_f
\end{equation}
in probability, where the {asymptotic mean vector} $a_f$ equals
\[
a_f=Ef \bigl(X_{t_1}^{st}, \ldots,
X_{t_r}^{st} \bigr).
\]
If, in addition, the initial distribution is such that for
some positive $\varepsilon$
%
\begin{equation}
\label{mu} \int_\XX \bigl(x^{-\gamma_i-\eps}+x^{\delta_i+\eps}
\bigr)\mu (\mathrm{d}x)<\infty ,\qquad i=1, \ldots, k,
\end{equation}
then (\ref{LLNd}) holds true in the mean sense.

2. (CLT). Assume in addition that there exists $%
\varepsilon>0$ such that
%
\begin{equation}
\label{squareint} E\bigl\llVert f \bigl(X_{t_1}^{st}, \ldots,
X_{t_r}^{st} \bigr)\bigr\rrVert^{2+%
\varepsilon}<\infty.
\end{equation}
Then
%
\begin{equation}
\label{CLTd} {\frac{1}{\sqrt n}}\sum_{l=1}^n
\bigl(f (X_{t_1+l}, \ldots, X_{t_r+l} )%
-a_f \bigr)\Rightarrow\mathcal{N}\bigl(0, \Sigma_{f}^d
\bigr),
\end{equation}
where the components of the {asymptotic covariance matrix} $\Sigma_{f}^d$
equal
\[
\bigl(\Sigma_f^d\bigr)_{i,j}=\sum
_{l=-\infty}^\infty\Cov \bigl(f_i
\bigl(%
X_{t_1+l}^{st}, \ldots, X_{t_r+l}^{st}
\bigr), f_j \bigl(X_{t_1}^{st}, \ldots,
X_{t_r}^{st} \bigr)\bigr) ,\qquad i,j =1, \ldots, k.
\]
\end{theo}

\begin{theo}[(Continuous-time case)]\label{t34}  Let the components of a vector-valued
function $f\dvtx \mathbb{X}^r\to\Re^k$ satisfy (\ref{fi}) with $\gamma_i,
\delta_i$
satisfying (\ref{gammadelta1}) for every $i=1,\ldots, k$.

Then the following statements hold true.

1. (LLN). For arbitrary initial distribution $\mu$ of $X$,
%
\begin{equation}
\label{LLNc} {\frac{1}{T}}\int_0^Tf
(X_{t_1+t}, \ldots, X_{t_r+t} )\,\mathrm{d}t\to a_f
\end{equation}
in probability. If, in addition, the initial distribution is such that for
some positive $\varepsilon$
%
\begin{equation}
\label{mu1} \int_\XX \bigl(x^{-(\gamma_i-1)\vee
0-\varepsilon}+x^{\delta_i+\varepsilon}
\bigr)\mu(\mathrm{d}x)<\infty,\qquad i=1,\ldots, k,
\end{equation}
then (\ref{LLNc}) holds true in the mean sense.

2. (CLT). Assume in addition that
%
\begin{equation}
\label{gammadelta2} \gamma_i< {\frac{\alpha}{4}}+{
\frac{1}{2}},\qquad \delta_i< {\frac
{\beta
}{4}},\qquad i=1,
\ldots, k.
\end{equation}
Then, for arbitrary initial distribution $\mu$ of $X$,
%
\begin{equation}
\label{CLTc} {\frac{1}{\sqrt T}}\int_0^T
\bigl(f (X_{t_1+t}, \ldots, X_{t_r+t} )-a_f%
\bigr)\,\mathrm{d}t\Rightarrow\mathcal{N}\bigl(0, \Sigma_{f}^c
\bigr),
\end{equation}
where the components of the {asymptotic covariance matrix} $\Sigma_{f}^c$
equal
%
\begin{equation}
\label{sigmaf} \bigl(\Sigma_f^c
\bigr)_{i,j}=\int_{-\infty}^\infty \Cov
\bigl(f_i \bigl(%
X_{t_1+t}^{st}, \ldots,
X_{t_r+t}^{st} \bigr), f_j \bigl(X_{t_1}^{st},
\ldots, X_{t_r}^{st} \bigr)\bigr) \,\mathrm{d}t,
\qquad i,j =1, \ldots, k.
\nonumber
\end{equation}
\end{theo}

For the limit theorems above, respective functional versions are
available, as well. In order to keep the exposition reasonably short,
we formulate here only one functional limit theorem of such a kind,
which corresponds to the CLT (\ref{CLTc}).

\begin{theo}
\label{tA} Let the components of a vector-valued function $f\dvtx \mathbb
{X}^r\to
\Re^k$ satisfy (\ref{fi}) with
%
\begin{equation}
\label{gammadelta3} \gamma_i< {\frac{\alpha}{4}},\qquad
\delta_i<{\frac{\beta}{4}},\qquad i=1,\ldots, k.
\end{equation}

Then
%
\begin{equation}
\label{CLTfunc} Y_T(\cdot)\equiv{\frac{1}{\sqrt T}}\int
_0^{ T} \bigl(f (X_{t_1+t}, \ldots,
X_{t_r+t} )-a_f \bigr)\,\mathrm{d}t\Rightarrow B,\qquad T
\to\infty
\end{equation}
weakly in $C([0,1])$, where $B$ is the Brownian motion in $\Re^k$ with the
covariance matrix of $B(1)$ equal to $\Sigma_{f}^c$.
\end{theo}

\section{Examples and statistical applications}

\subsection{Examples}

In this section, we illustrate the above limit
theorems and use them to derive the asymptotic properties of \emph
{empirical mixed moments}
\[
\overline{m}_{\upsilon,\chi, c} (t) = \frac{1}{T} \int_0^T
X_{s}^{\upsilon} X_{t+s}^{\chi} \,\mathrm{d}s,
\qquad \overline{m}_{\upsilon,\chi
,d} (t) = \frac{1%
}{n} \sum
_{l=1}^n X_{l}^{\upsilon}X_{t+l}^{\chi},
\qquad t>0
\]
both in the continuous-time and in the discrete-time settings. Below we
use statistical terminology because such functionals are particularly
important for the statistic inference. For instance, usual \emph
{empirical moments}
%
\begin{equation}
\label{empmom} \overline{m}_{\upsilon,c} = \frac{1}{T} \int
_0^T X_{s}^{\upsilon} \,
\mathrm{d}s,\qquad \overline{m}_{\upsilon,d} = \frac{1%
}{n} \sum
_{l=1}^n X_{l}^{\upsilon}
\end{equation}
equal the empirical mixed moments with $\chi=0$, and \emph{empirical
covariances}
%
\begin{eqnarray}
\label{empcov} \overline{R}_{c}(t)&=&\frac{1}{T} \int
_0^T X_sX_{t+s}\,
\mathrm{d}s- \biggl(\frac{1}{T} \int_0^T
X_{s}\,\mathrm{d}s \biggr)^2,
\nonumber\\[-8pt]\\[-8pt]
\overline{R}_{d}(t)&=&\frac{1}{n} \sum
_{l=1}^n X_{l}X_{t+l}-
\Biggl(\frac{1%
}{n} \sum_{l=1}^n
X_{l} \Biggr)^2,
\nonumber
\end{eqnarray}
can be written as
%
\begin{equation}
\label{empcov2} \overline{R}_{c}(t)=\overline{m}_{1,1,c}(t)-
(\overline {m}_{1,c} )^2,\qquad \overline{R}_{d}(t)=
\overline{m}_{1,1,d}(t)- (\overline {m}_{1,d} )^2.
\end{equation}

Denote $\upsilon_-=-(\upsilon\wedge0), \upsilon_+=\upsilon\vee0.$

\begin{ex}[(Discrete-time case)]\label{e31} Let there exist
$p,q>1$ with $1/p+1/q=1$ such that
%
\begin{equation}
\label{u1} \{p\upsilon, q\chi\}\subset \biggl(-{\frac{\alpha}{2}}+1, {
\frac
{\beta
}{2}} \biggr).
\end{equation}
Then for arbitrary initial distribution $\mu$
of $X$ the discrete-time empirical mixed moment $\overline
{m}_{\upsilon
,\chi,d}(t)$ is a $P$-consistent estimator of the parameter
\[
m_{\upsilon,\chi}(t)=E \bigl(X_{0}^{st}
\bigr)^\upsilon \bigl(X_{t}^{st} \bigr)^\chi.
\]

If, in addition, the initial distribution $\mu$ satisfies
\[
\int_0^1 x^{-(p\upsilon_-)\vee(q\chi_-)-\varepsilon}\mu(\mathrm{d}x)+
\int_1^\infty x^{(p\upsilon_+)\vee(q\chi_+)+\varepsilon}\mu(\mathrm{d}x)<
\infty
\]
for some $\eps>0$,
then $\overline{m}_{\upsilon,\chi,d}(t)$ is an asymptotically unbiased
estimator of $m_{\upsilon,\chi}(t)$.\eject

Under the assumption
%
\begin{equation}
\label{u2} \{p\upsilon, q\chi\}\subset \biggl(- \biggl({\frac{\alpha}{2}}%
-1 \biggr)\wedge \biggl({\frac{\alpha}{4}} \biggr), {\frac{\beta
}{4}}
\biggr)
\end{equation}
for arbitrary initial distribution $\mu$ of $X$ the discrete-time
empirical mixed moment $\overline{m}_{\upsilon,\chi,d}(t)$ is an
asymptotically normal
estimator of $m_{\upsilon,\chi}(t)$; that is,
\[
\sqrt n \bigl(\overline{m}_{\upsilon,\chi,d}(t)-m_{\upsilon,\chi
}(t) \bigr)
\Rightarrow\mathcal{N}\bigl(0, \sigma_{\upsilon,\chi,d}^2(t)\bigr),
\qquad n\to \infty
\]
with
\[
\sigma_{\upsilon,\chi,d}^2(t)=\sum_{l=-\infty}^\infty
 \Cov \bigl( \bigl(X_{l}^{st} \bigr)^\upsilon
\bigl(X_{t+l}^{st} \bigr)^\chi,
\bigl(X_{0}^{st} \bigr)^\upsilon
\bigl(X_{t}^{st} \bigr)^\chi \bigr).
\]
\end{ex}

These results follow immediately from Theorem \ref{t32} with $k=1,$
$r=2$, and
\[
f(x_1, x_2)=x_1^\upsilon
x_2^\chi.
\]
Indeed, by the Young inequality,
\[
f(x_1,x_2)\leq\frac{x_1^{p\upsilon}}{ p}+ \frac{x_2^{q\chi}}{ q}.
\]
Then (\ref{fi}) holds true with $\gamma=(p\upsilon_-)\vee(q\chi_-)$
and $\delta=(p\upsilon_+)\vee(q\chi_+)$. Respectively,
(\ref{u1}) coincides with the assumption (\ref{gammadelta}), imposed on
$\gamma, \delta$ in Theorem \ref{t32}. The additional integrability
assumption (\ref{squareint}) now is equivalent to the following: for
some positive $\eps$,
\[
-2(p\upsilon_-)\vee(q\chi_-)-\eps>-\frac{\alpha}{2},\qquad 2(p\upsilon_+)
\vee(q\chi_+)+\eps<\frac{\beta}{2}.
\]
Clearly, this means that $\{p\upsilon, q\chi\}\subset(-\alpha/4,
\beta
/4)$, which together with (\ref{u1}) gives (\ref{u2}).

Similarly, using Theorem \ref{t34} under the same choice of $f,\gamma,
\delta$ we obtain the following.

\begin{ex}[(Continuous-time case)]\label{e32}
Let there exist $p,q>1$ with $1/p+1/q=1$ such that
%
\begin{equation}
\label{u3} \{p\upsilon, q\chi\}\subset \biggl(-{\frac{\alpha}{2}}, {
\frac
{\beta
}{2}} \biggr).
\end{equation}
Then for arbitrary initial distribution $\mu$
of $X$ the continuous-time empirical mixed moment $\overline
{m}_{\upsilon,\chi,c}(t)$ is a $P$-consistent estimator of the $m_{\upsilon,\chi}(t)$.

If, in addition, the initial distribution $\mu$ satisfies
\[
\int_0^1 x^{-((p\upsilon_-)\vee(q\chi_-)-1)_+-\varepsilon}\mu (\mathrm{d}x)+
\int_1^\infty x^{(p\upsilon_+)\vee(q\chi_+)+\varepsilon}\mu(\mathrm{d}x)<
\infty
\]
for some $\eps>0$,
then $\overline{m}_{\upsilon,\chi,c}(t)$ is an asymptotically unbiased
estimator of $m_{\upsilon,\chi}(t)$.

Under the assumption
%
\begin{equation}
\label{u4} \{p\upsilon, q\chi\}\subset \biggl(-{\frac{\alpha}{4}}-{
\frac
{1}{2}}, {%
\frac{\beta}{4}} \biggr)
\end{equation}
for arbitrary initial distribution $\mu$ of $X$ the continuous-time
empirical mixed moment $\overline{m}_{\upsilon,\chi,c}(t)$ is an
asymptotically normal
estimator of $m_{\upsilon,\chi}(t)$; that is,
\[
\sqrt T \bigl(\overline{m}_{\upsilon,\chi,c}(t)-m_{\upsilon,\chi
}(t) \bigr)
\Rightarrow\mathcal{N}\bigl(0, \sigma_{\upsilon,\chi,c}^2(t)\bigr),
\qquad T\to \infty
\]
with
\[
\sigma_{\upsilon,\chi,c}^2(t)=\int_{-\infty}^\infty
\Cov \bigl( \bigl(X_{s}^{st} \bigr)^\upsilon
\bigl(X_{t+s}^{st} \bigr)^\chi,
\bigl(X_{0}^{st} \bigr)^\upsilon
\bigl(X_{t}^{st} \bigr)^\chi \bigr)\,\mathrm{d}s.
\]
\end{ex}

The following statements can be obtained easily either by taking in the
above examples $\chi=0$ and $p>1$ close enough to $1$, or by using
Theorem \ref{t32} and Theorem \ref{t34} with $k=r=1, f(x)=x^\upsilon$,
and $\gamma=\upsilon_-, \delta=\upsilon_+$.

\begin{ex}[(Empirical moments)]\label{e33} The discrete-time
empirical moment $\overline{m}_{\upsilon,d}$, considered as an
estimator of the parameter
\[
m_\upsilon=E \bigl(X_{0}^{st} \bigr)^\upsilon=
\int_{\XX}x^\upsilon\pi(\mathrm{d}x),
\]
has the following properties:
\begin{enumerate}[(iii)]
\item[(i)] if
%
\begin{equation}
\label{u5} \upsilon\in \biggl(-{\frac{\alpha}{2}}+1, {\frac{\beta}{2}}
\biggr),
\end{equation}
then $\overline{m}_{\upsilon,d}$ is $P$-consistent;
\item[(ii)] if, in addition, the initial distribution $\mu$ satisfies
\[
\int_0^1 x^{-\upsilon_--\varepsilon}\mu(\mathrm{d}x)+
\int_1^\infty x^{\upsilon
_++\varepsilon}\mu(\mathrm{d}x)<
\infty
\]
for some $\eps>0$,
then $\overline{m}_{\upsilon,d}$ is asymptotically unbiased;

\item[(iii)] if
%
\begin{equation}
\label{u6} \upsilon\in \biggl(- \biggl({\frac{\alpha}{2}}%
-1 \biggr)
\wedge \biggl({\frac{\alpha}{4}} \biggr), {\frac{\beta
}{4}} \biggr),
\end{equation}
then $\overline{m}_{\upsilon,d}$ is asymptotically normal.
\end{enumerate}

Similarly, the continuous-time empirical moment $\overline
{m}_{\upsilon
,c}$, considered as an estimator of the same parameter, satisfies the following:
\begin{enumerate}[(iii)]
\item[(i)] if
%
\begin{equation}
\label{u7} \upsilon\in \biggl(-{\frac{\alpha}{2}}, {\frac{\beta}{2}}
\biggr),
\end{equation}
then $\overline{m}_{\upsilon,c}$ is $P$-consistent;
\item[(ii)] if, in addition, the initial distribution $\mu$ satisfies

\[
\int_0^1 x^{-(\upsilon_--1)_+-\varepsilon}\mu(\mathrm{d}x)+
\int_1^\infty x^{\upsilon_+
+\varepsilon}\mu(\mathrm{d}x)<
\infty
\]
for some $\eps>0$,
then $\overline{m}_{\upsilon,d}$ is asymptotically unbiased;

\item[(iii)] if
%
\begin{equation}
\label{u8} \upsilon\in \biggl(-{\frac{\alpha}{4}}-{\frac{1}{2}},
{%
\frac{\beta}{4}} \biggr),
\end{equation}
then $\overline{m}_{\upsilon,c}$ is asymptotically normal.
\end{enumerate}
\end{ex}

Comparing (\ref{u5}) with (\ref{u7}) and (\ref{u6}) with (\ref{u8}),
one can see clearly the difference
between the conditions of Theorem \ref{t34} and
the conditions of
Theorem \ref{t32}. The particularly interesting case here is
\[
\upsilon\in \biggl(-{\frac{\alpha}{4}}-{\frac{1}{2}}, -{\frac
{\alpha
}{4}}
\biggr].
\]
In this case, the function $f(x)=x^{\upsilon}$ satisfies conditions of
Theorem \ref{t34} with $r=k=1$, while the additional integrability
assumption (\ref{squareint}) in Theorem \ref{t32} fails because $f$ is
not square
integrable w.r.t. $\pi$. This observation reveals a new effect, already
mentioned in the
Introduction, which seemingly has not been observed in the literature
before: a functional $f$, which is not square integrable w.r.t. the
invariant distribution, still may lead to the CLT in its continuous-time
form~(\ref{CLTc}).

\begin{ex}[(Empirical covariances)]\label{e34}  Both the
discrete-time empirical covariance $\overline{R}_{d}(t)$ and the
continuous-time empirical covariance $\overline{R}_{c}(t)$, considered
as estimators of the parameter
\[
R(t)=\Cov \bigl(X_{t}^{st}, X_{0}^{st}
\bigr),
\]
have the following properties:
\begin{enumerate}[(iii)]
\item[(i)] if $\beta>4$ then $\overline{R}_{d}(t)$ and $\overline
{R}_{c}(t)$ are
$P$-consistent;
\item[(ii)] if, in addition, the initial distribution $\mu$ satisfies
\[
\int_1^\infty x^{2+\varepsilon}\mu(\mathrm{d}x)<
\infty
\]
for some $\eps>0$,
then $\overline{R}_{d}(t)$ and $\overline{R}_{c}(t)$ are asymptotically
unbiased;

\item[(iii)]
if $\beta>8$ then $\overline{R}_{d}(t)$ and $\overline{R}_{c}(t)$ are
asymptotically normal.
\end{enumerate}
\end{ex}

These results follow from the representation (\ref{empcov2}) and
Theorems \ref{t32}, \ref{t34} with $k=r=2$, $f=(f_1,f_2)$,
\[
f_1(x_1,x_2)=x_1,\qquad
f_2(x_1,x_2)=x_1x_2.
\]
Similarly to Example \ref{e31} and Example \ref{e32} (in this
particular case one should take $p=q=2$), one can verify that both
$(\overline{m}_{1,d}, \overline{m}_{1,1,d}(t))$ and $(\overline
{m}_{1,c}, \overline{m}_{1,1,c}(t))$ are $P$-consistent if $\beta>4$
and asymptotically normal if $\beta>8$, when considered as estimators
of the vector parameter $({m}_{1},{m}_{1,1}(t))$. Then properties (i)
and (iii) follow by the continuity mapping theorem and the functional
delta method (see \cite{Serfling}, Theorem 3.3.A). Under the additional
integrability assumption on $\mu$ both $\overline{m}_{1,1,d}(t)$ and
$\overline{m}_{1,1,c}(t)$ are asymptotically unbiased. On the other
hand, under the same assumption both $(\overline{m}_{1,d})^2$ and
$(\overline{m}_{1,c})^2$ are uniformly integrable w.r.t. $P_\mu$; this
follows from the H\"older inequality and Corollary \ref{c31}:
\[
E_\mu(\overline{m}_{1,d})^{2+\eps}=E_\mu
\Biggl(\frac{1%
}{n} \sum_{l=1}^n
X_{l} \Biggr)^{2+\eps}\leq\frac{1%
}{n} \sum
_{l=1}^n E_\mu X_{l}^{2+\eps}
\leq C,
\]
the inequality for the continuous-time case is similar and omitted.
This implies that $(\overline{m}_{1,d})^2$ and $(\overline{m}_{1,c})^2$
are asymptotically unbiased, which completes the proof of the property (ii).

Similarly, the properties of the empirical estimates of the
vector-valued parameters of the type
$({m}_{\upsilon_1}, \ldots, {m}_{\upsilon_k})$ or
$({m}_{\upsilon_1}, \ldots, {m}_{\upsilon_k}, R(t))$ can be derived. For
such parameters, the component-wise properties of $P$-consistency and
asymptotic unbiasedness are already studied in the previous examples.
Hence, in the following example, we address the asymptotic normality only.

\begin{ex}[(Multivariate estimators)]\label{e35} \emph{I. (Discrete-time
case).} Let
\[
\upsilon_1, \ldots, \upsilon_k\in \biggl(- \biggl({
\frac{\alpha}{2}}%
-1 \biggr)\wedge \biggl({\frac{\alpha}{4}}
\biggr), {\frac{\beta
}{4}} \biggr).
\]
Then, for arbitrary initial distribution $\mu$ of $X$, the estimator
$\overline{m}_{\upsilon_1, \ldots, \upsilon_k,d}=(\overline
{m}_{\upsilon_1,d},
\ldots,\break \overline{m}_{\upsilon_k,d})$ of the vector-valued parameter $
m_{\upsilon_1, \ldots, \upsilon_k}=({m}_{\upsilon_1}, \ldots,
{m}_{\upsilon_k})
$ is asymptotically normal; that is,
\[
\sqrt n (\overline{m}_{\upsilon_1, \ldots, \upsilon
_k,d}-{m}_{\upsilon_1,
\ldots, \upsilon_k} )\Rightarrow
\mathcal{N}(0, \Sigma),\qquad n\to \infty
\]
with some positive semi-definite matrix $\Sigma$.

If, in addition, $\beta>8$%
, then $(\overline{m}_{\upsilon_1,d}, \ldots, \overline{m}_{\upsilon
_k,d},\overline{R}_{d}(t))$ is an asymptotically normal estimator of
$({m}%
_{\upsilon_1}, \ldots, {m}_{\upsilon_k}, R(t))$ for any $t>0$.

\emph{II. (Continuous-time case).} Let
\[
\upsilon_1, \ldots, \upsilon_k\in \biggl(-{
\frac{\alpha}{4}}-{\frac{1}{2}}, {%
\frac{\beta}{4}}
\biggr).
\]
Then, for arbitrary initial distribution $\mu$ of $X$, the estimator $%
\overline{m}_{\upsilon_{1}, \ldots, \upsilon_k,c}=(\overline{m}%
_{\upsilon_1,c}, \ldots,\break \overline{m}_{\upsilon_k,c})$ of the
vector-valued\vadjust{\goodbreak}
parameter $m_{\upsilon_1, \ldots, \upsilon_k}=({m}_{\upsilon_1},
\ldots,
{m}_{\upsilon_k})$ is asymptotically normal; that is,
\[
\sqrt T (\overline{m}_{\upsilon_1, \ldots, \upsilon
_k,c}-{m}_{\upsilon_1,
\ldots, \upsilon_k} )\Rightarrow
\mathcal{N}(0, \Sigma),\qquad T\to \infty
\]
with some positive semi-definite matrix $\Sigma$.

If, in addition, $\beta>8$, then $(\overline{m}_{\upsilon_1,c},
\ldots,
\overline{m}_{\upsilon_k,c}, \overline{R}_{c}(t))$ is an asymptotically
normal estimator of $({m}_{\upsilon_1}, \ldots, {m}_{\upsilon_k},
R(t))$ for
any $t>0$.
\end{ex}

\subsection{Parameter estimation for the Fisher--Snedecor diffusion}\label{s34}

In this section, we give an application of the above
results to
the parameter estimation of the Fisher--Snedecor diffusion. We use the
method of moments and the asymptotic properties of
the {empirical moments} (\ref{empmom}) and the {empirical covariances}
(\ref{empcov}), exposed in Examples \ref{e33}--\ref{e35},
in order to
provide the statistical analysis of the \emph{autocorrelation
parameter} $%
\theta$ and the \emph{shape parameters} $\alpha,$ $\beta$, and
$\kappa$ of
the Fisher--Snedecor diffusion. We put
%
\begin{eqnarray}
\label{estt} %
\widehat{\alpha}_c&=&\frac{2(
\overline{m}_{-1,c}\overline{m}_{1,c}\overline{m}_{2,c}-%
\overline{m}_{1,c}^2)}{
\overline{m}_{-1,c}\overline{m}_{1,c}\overline{m}_{2,c}-2\overline
{m}_{2,c}+%
\overline{m}_{1,c}^2 },
\qquad \widehat{\beta}_c=\frac{4\overline{m}_{-1,c}(\overline{m}_{2,c}-%
\overline{m}_{1,c}^2)}{
\overline{m}_{-1,c}\overline{m}_{2,c}-2\overline{m}_{-1,c}%
\overline{m}_{1,c}^2+ \overline{m}_{1,c}},\quad
\nonumber
\\[-8pt]
\\[-8pt]
\widehat{\kappa}_c&=&\frac{4\overline{m}_{-1,c}\overline{m}_{1,c}(%
\overline{m}_{2,c}-\overline{m}_{1,c}^2)}{
\overline{m}_{-1,c}\overline{m}_{2,c}-2\overline{m}_{-1,c}%
\overline{m}_{1,c}^2+ \overline{m}_{1,c}},\qquad\widehat{\theta
}_c=-\frac{1}{
t}\log \biggl(\frac{\overline{R}_c(t)}{
\overline{m}_{2,c}-\overline{m}_{1,c}^2} \biggr)
\nonumber
\end{eqnarray}
for a given $t>0$, and define $\widehat{\alpha}_d, \widehat{\beta}_d,
\widehat{\kappa}_d, \widehat{\theta}_d$ by similar relations with
$\overline{%
m}_{i,d}, i=-1, 1,2$, and $\overline{R}_d(t)$ instead of $\overline{m}
_{i,c}, i=-1, 1,2$, and $\overline{R}_c(t)$, respectively.

\begin{theo}
\label{t36} Let $\beta>8$. Then, for arbitrary initial distribution
of the
Fisher--Snedecor diffusion, $(\widehat{\alpha}_c, \widehat{\beta}_c,
\widehat{%
\kappa}_c, \widehat{\theta}_c)$ is a $P$-consistent and asymptotically
normal estimator of the parameter $(\alpha, \beta, \kappa, \theta
)$; that
is,
\[
\sqrt{T}(\widehat{\alpha}_c-\alpha, \widehat{\beta}_c-
\beta, \widehat {\kappa}%
_c-\kappa, \widehat{
\theta}_c-\theta)\Rightarrow\mathcal{N}%
\bigl(0,
\Sigma_c(\alpha, \beta, \kappa, \theta)\bigr),\qquad T\to\infty.
\]

For the estimator $(\widehat{\alpha}_d, \widehat{\beta}_d, \widehat
{\kappa}%
_d, \widehat{\theta}_d),$ the similar statement holds true under the
additional assumption $\alpha>4$. In that case,
\[
\sqrt{n}(\widehat{\alpha}_d-\alpha, \widehat{\beta}_d-
\beta, \widehat {\kappa}%
_d-\kappa, \widehat{
\theta}_d-\theta)\Rightarrow\mathcal{N}%
\bigl(0,
\Sigma_d(\alpha, \beta, \kappa, \theta)\bigr),\qquad n\to\infty.
\]
\end{theo}

The matrices $\Sigma_c(\alpha, \beta, \kappa, \theta)$, $\Sigma_d(\alpha,
\beta, \kappa, \theta)$ are completely identifiable. To keep the current
paper reasonably short, we postpone their explicit calculation, together
with a more detailed discussion of the statistical aspects, to the
subsequent paper \cite{KulikLeonenko}.

\begin{rem}
The estimators (\ref{estt}) can be simplified significantly if either exact
values of some parameters $\alpha, \beta, \kappa$ are known, or these
parameters possess some functional relation. Let, for instance, $\kappa
={%
\beta/(\beta-2)}$; this particular case is of a separate interest because
the invariant distribution $\pi$ then coincides with the Fisher--Snedecor
distribution $\mathcal{F S}(\alpha, \beta)$. In this case, one can replace
in (\ref{estt}) the identities for $\widehat{\alpha}_c, \widehat
{\beta}_c$
by either
%
\begin{equation}
\label{estd} \widehat{\alpha}_c = \frac{2 \overline{m}_{1,c}^{ 2}}{\overline{m}%
_{2,c}(2 - \overline{m}_{1,c}) - \overline{m}_{1,c}^{ 2}},\qquad
\widehat{%
\beta}_c = \frac{2\overline{m}_{1,c}}{\overline{m}_{1,c}-1}
\end{equation}
or
%
\begin{equation}
\label{estc} \widehat{\alpha}_c = \frac{2\overline{m}_{-1,c}}{\overline{m}_{-1,c}-1},\qquad
\widehat{\beta}_c = \frac{2\overline{m}_{1,c}}{\overline{m}_{1,c}-1}.
\end{equation}
For the estimator $(\widehat{\alpha}_c, \widehat{\beta}_c, \widehat
{\theta}%
_c),$ defined in such a way, and its discrete-time analogue $(\widehat
{\alpha%
}_d, \widehat{\beta}_d, \widehat{\theta}_d)$, the statements of Theorem
\ref{t36} hold true; see more detailed discussion in~\cite{KulikLeonenko}.
\end{rem}

\section{Proofs}\label{s4}

\subsection{\texorpdfstring{Proof of Proposition \protect\ref{p21}}{Proof of Proposition 3.1}}

\textit{Statement 1.} Let the initial value $X_0=x\in\mathbb{X}$ be fixed.
Note that the process
%
\begin{equation}
\label{Hf} H^{\phi,X}_t=\phi(X_t)-\int
_0^t \mathcal{A} \phi(X_s)\,
\mathrm{d}s,\qquad t\in \Re^+,
\end{equation}
with $\mathcal{A} \phi$ defined by \eqref{Af}, is an $\mathbb{F}^X$-local
martingale w.r.t. the measure $P_x.$ The argument here is quite
standard, we
explain it briefly in order to keep the exposition self-sufficient.
Introduce the sequence of $\mathbb{F}^X$-stopping times $%
T_n=\inf\{t\dvt X_t\leq1/n\}, n\in\mathbb{N}$, and consider auxiliary
functions $%
\phi_n\in C^2(\Re)$ such that $\phi_n=\phi$ on $[1/n,\infty)$. For any
given $n\in
\mathbb{N}$, by the Ito formula (e.g., \cite{IkedaWat}, Chapter II, Theorem
5.1) we have that the process $H^{\phi_n,X}$, defined by the relation
\eqref{Hf}
with $\phi_n$ instead of $\phi$, is an $\mathbb{F}^X$-local
martingale. This
means that, for any given $n\in\mathbb{N}$, there exists a sequence of
$%
\mathbb{F}^X$-stopping times $T_{n,m}, m\in\mathbb{N}$ such that every
process
\[
t\mapsto H^{\phi_n,X}(t\wedge T_{n,m}),\qquad m\in\mathbb{N}
\]
is an $\mathbb{F}^X$-martingale w.r.t. the measure $P_x,$ and
\[
T_{n,m}\to\infty,\qquad m\to\infty, P_x\mbox{-a.s.}
\]
The last relation provides that for every $n\in\mathbb{N}$ there
exists $m_n
$ such that
\[
P_x(T_{n,m_n}\leq n)<2^{-n}.
\]
Consequently, by the Borel--Cantelli lemma,
\[
T_{n,m_n}\to\infty,\qquad n\to\infty, P_x\mbox{-a.s.}
\]
On the other hand, since the point $0$ is unattainable for $X$, we have
$%
T_n\to\infty$ $P_x$-a.s. Consequently, for $S_n=T_n\wedge T_{n,m_n},
n\in
\mathbb{N}$ we have
\[
S_{n}\to\infty,\qquad n\to\infty, P_x\mbox{-a.s.}
\]
By the Doob optional sampling theorem, the process
\[
t\mapsto H^{\phi_n,X}(t\wedge S_{n})
\]
is an $\mathbb{F}^X$-martingale w.r.t. the measure $P_x$. On the other
hand, the processes $H^{\phi_n,X}$ and $H^{\phi,X}$ coincide up to the
time moment
$T_n$ because the values of $\phi_n$ and its derivatives on $[1/n,
\infty)$
coincide with respective values of $\phi$. Hence, the process
\[
t\mapsto H^{\phi,X}(t\wedge S_{n})
\]
is an $\mathbb{F}^X$-martingale w.r.t. the measure $P_x$, which completes
the proof of the fact that $H^{\phi,X}$ is a $\mathbb{F}^X$-local
martingale.

Next, we show that the function $\mathcal{A} \phi$ defined by (\ref{Af})
satisfies (\ref{lyap}) for properly chosen positive $u,v,c,C$. We have
for $x
$ large enough:
%
\begin{eqnarray}
\label{xlarge} %
\mathcal{A} \phi(x)&=& -\theta\delta(x-\kappa)
x^{\delta-1}+\theta\delta(\delta-1)x \biggl(\frac{x}{\beta
/2-1}+
\frac{\kappa
}{
\alpha/2} \biggr)x^{\delta-2}
\nonumber
\\[-8pt]
\\[-8pt]
&=&-\theta\delta \phi(x) \biggl[ \biggl(1-\frac{\kappa}{ x} \biggr) -(\delta-1)
\biggl(\frac{1}{
\beta/2-1}+\frac{\kappa}{ x\alpha/2} \biggr) \biggr].
\nonumber
\end{eqnarray}
The term $ [\cdots ]$ tends to $1-{\frac{\delta-1}{\beta/2-1}}$
as $%
x\to\infty$, and it was assumed that $\delta<\beta/2$. Hence (\ref{lyap})
holds true for any $x> v$ assuming $v>0$ is chosen large enough and
$c>0$ is
chosen small enough.

We have for $x$ small enough:
%
\begin{eqnarray}
\label{xsmall} %
\mathcal{A} \phi(x)&=& \theta\gamma(x-\kappa)
x^{-\gamma-1}+\theta\gamma(\gamma+1)x \biggl(\frac{x}{\beta
/2-1}+
\frac{\kappa
}{
\alpha/2} \biggr)x^{-\gamma-2}
\nonumber
\\[-8pt]
\\[-8pt]
&=&-\theta\gamma\phi(x) \biggl\{ \biggl(\frac{\kappa}{ x}-1 \biggr)- (\gamma+1)
\biggl(\frac{1}{
\beta/2-1}+\frac{\kappa}{ x\alpha/2} \biggr) \biggr\}.
\nonumber
\end{eqnarray}
The term $ \{\cdots \}$ is equivalent to
\[
{\frac{\kappa}{x}} \biggl(1-{\frac{\gamma+1}{\alpha/2}} \biggr)
\]
as $x\to0+$, and it tends to $+\infty$ because it was assumed that $%
\gamma+1<\alpha/2$. Hence, (\ref{lyap}) holds true for any $x\in(0,u)$
assuming $u, c>0$ are chosen small enough. Finally, for given $u,v,c$
(\ref{lyap}) holds true for $x\in[u,v]$ under appropriate choice of (large) $C$.

Finally, we show that the process (\ref{Hf}) is an $\mathbb{F}^X$-martingale. This proof is quite standard, again. For any $n\in\mathbb{N}$,
we have
%
\begin{equation}
\label{11} E_x H^{\phi,X}(t\wedge S_{n})=\phi(x),
\qquad t\geq0;
\end{equation}
here $S_n, n\in\mathbb{N}$ is the sequence of stopping times
constructed in
the first part of the proof. Recall that it is supposed that $\phi
(x)\geq1$, and therefore $\phi(x)$ is positive.
This, together with (\ref{lyap}), provides that $[\mathcal{A} \phi
]_+(x)=(%
\mathcal{A} \phi(x))\vee0$ is a bounded function. Then
\[
E_x \phi(X_{t\wedge S_{n}})=\phi(x)+E_x\int
_0^{t\wedge
S_{n}}\mathcal{A} \phi(X_s)\,
\mathrm{d}s\leq\phi(x)+t\sup_{x^{\prime}}[\mathcal{A} \phi ]_+
\bigl(x^{\prime}\bigr),\qquad t\geq0, n\in\mathbb{N}.
\]
Consequently, we have from (\ref{11}) that for any $T\geq0$
%
\begin{equation}
\label{12} \sup_{t\leq T}\sup_{n\in\mathbb{N}} E_x
\phi(X_{t\wedge
S_{n}})<\infty.
\end{equation}
Denote $[%
\mathcal{A} \phi]_-(x)=(-\mathcal{A} \phi(x))\vee0$; then (\ref{11})
can be
written as
\[
E_x\int_0^{t\wedge S_{n}}[\mathcal{A}
\phi]_-(X_s)\,\mathrm{d}s=\phi(x)-E_x \phi
(X_{t\wedge
S_{n}})+E_x\int_0^{t\wedge S_{n}}[
\mathcal{A} \phi]_+(X_s)\,\mathrm{d}s.
\]
Combined with (\ref{12}) and the fact that $[\mathcal{A} \phi]_+$ is bounded,
this yields
\[
E_x\int_0^{t}[\mathcal{A}
\phi]_-(X_s)\,\mathrm{d}s<\infty.
\]
In particular, the Lebesgue dominated convergence theorem and
boundedness of
$[\mathcal{A} \phi]_+$ provide that the sequence
\[
\int_0^{t\wedge S_{n}}\mathcal{A} \phi(X_s)
\,\mathrm{d}s,\qquad n\in\mathbb{N}
\]
is uniformly integrable w.r.t. $P_x$.

Note that the above argument can be repeated with the function $\phi$ replaced
by the function $\tilde\phi=\phi^\upsilon$, where $%
\upsilon>1$ is chosen in such a way that
\[
\upsilon\gamma<{\frac{\alpha}{2}}-1,\qquad \upsilon\delta<{
\frac{\beta}{2}}.
\]
Then, similarly to (\ref{12}), we will have
%
\begin{equation}
\label{13} \sup_{t\leq T}\sup_{n\in\mathbb{N}} E_x \bigl(
\phi(X_{t\wedge
\tilde
S_{n}}) \bigr)^\upsilon<\infty
\end{equation}
with some sequence of stopping times $\tilde S_n$ such that $\tilde
S_n\to
\infty$ $P_x$-a.s. This means that the sequence $\phi(X_{t\wedge
S_{n}\wedge\tilde S_n}), n\in\mathbb{N}$ of
the processes on $[0,T]$ is uniformly integrable w.r.t. $P_x$, and
hence the
sequence $H^{\phi,X}(t\wedge S_{n}\wedge\tilde S_n), n\in\mathbb
{N}$ is
uniformly integrable, as well. Then $H^{\phi,X}$ is a martingale as an a.s.
limit of a uniformly integrable sequence of martingales.

\textit{Statement 2.} Take a segment $[w,z]\in\mathbb{X}$ such that $
[u,v]\subset(w,z)$, and consider the process $X^{[w,z]}$ obtained from $X$
by killing at the exit from $(w,z)$. Clearly, for any $x$ inside
$(w,z)$ the
transition probability $P_t(x,\mathrm{d}y)$ is minorized by the transition
probability $P_t^{[w,z]}(x,\mathrm{d}y)$ of the process $X^{[w,z]}$. The latter
function is the fundamental solution to the Cauchy problem for the linear
2nd order parabolic equation
\[
\partial_t u (x,y)=\mathcal{L} u(t,x),\qquad x\in(w,z),
u(t,w)=u(t,z)=0, t> 0,
\]
where
\[
\mathcal{L}=a(x)\partial_x+{\tfrac{1}{2}}\sigma^2(x)
\partial^2_{xx}.
\]
Because the coefficients $a,\sigma$ are smooth in $[w,z]$ and $\sigma
$ is
positive, the general analytic results from the theory of linear 2-nd order
parabolic equations (e.g., \cite{LadSolUr}, Chapter IV, Sections
11--14) yield
representation
\[
P_t^{[w,z]}(x,\mathrm{d}y)=Z_t(x,y)\,\mathrm{d}y
\]
with a continuous function $Z\dvtx (0,+\infty)\times(w,z)\times(w,z)\to
[0,\infty)$. Because $Z$ is continuous and is not an identical zero, there
exist $t_1>0, x_1\in(w,z), y_1\in(w,z)$, and $\varepsilon>0$ such that
\[
c_1:=\inf_{|x-x_1|\leq\varepsilon, |y-y_1|\leq\varepsilon}Z_{t_1}(x,y)>0.
\]
In other words, we have constructed $t_1>0$ and segments $%
[u^{\prime},v^{\prime}]=[y_1-\varepsilon, y_1+\varepsilon]$ and $%
[u^{\prime\prime},v^{\prime\prime}]=[x_1-\varepsilon,
x_1+\varepsilon]$ such
that
%
\begin{equation}
\label{14} P_{t_1}(x,A)\geq P_{t_1}^{[w,z]}(x,A)\geq
c_1\int_A \mathrm{d}y
\end{equation}
for any $x\in[u^{\prime\prime},v^{\prime\prime}]$ and Borel measurable
set $%
A\subset[u^{\prime},v^{\prime}]$. Take $t_2>0$ and put $T=t_1+t_2$. The
Chapmen--Kolmogorov equation and (\ref{14}) yields for every $x\in
[u,v]$ and
Borel measurable set $A\subset[u^{\prime},v^{\prime}]$
\[
P_T(x,A)\geq \int_{[u^{\prime\prime},v^{\prime\prime}]}P_{t_1}
\bigl(x^{%
\prime},A\bigr)P_{t_2}\bigl(x,\mathrm{d}x^{\prime}
\bigr)\geq c_1 \inf_{x\in
[u,v]}P_{t_2}\bigl(x,
\bigl(u^{\prime\prime}, v^{\prime\prime}\bigr)\bigr)\int_{A}
\mathrm{d}y.
\]
The reason for us to replace in the last inequality the segment $%
[u^{\prime\prime},v^{\prime\prime}]$ by the open interval $%
(u^{\prime\prime},v^{\prime\prime})$ is that the indicator of this interval
can be obtained as a limit of an increasing sequence of continuous functions
$f_n\dvtx \mathbb{X}\to\Re^+, n\geq1$. The process $X$ is a Feller one; this
follows from the standard theorem on continuity of a solution to an SDE
w.r.t. its initial value, for example, \cite{GhSk}, Chapter II.
Therefore, every
function
\[
x\mapsto\int_\XX f_n(y)P_{t_2}(x,
\mathrm{d}y)
\]
is continuous, which implies that the function
\[
x\mapsto P_{t_2}\bigl(x,\bigl(u^{\prime\prime}, v^{\prime\prime}\bigr)
\bigr)
\]
is lower semicontinuous as a point-wise limit of an increasing sequence of
continuous functions. Then there exists $x_\diamondsuit\in[u,v]$ such that
\[
\inf_{x\in[u,v]}P_{t_2}\bigl(x,\bigl(u^{\prime\prime},
v^{\prime\prime}\bigr)\bigr)=P_{t_2}\bigl(x_\diamondsuit,u^{\prime\prime},
v^{\prime\prime}\bigr).
\]
On the other hand, for any $t>0, x\in\mathbb{X}$ the support of the measure
$P_t(x,\cdot)$ coincides with whole $\mathbb{X}$; because the diffusion
coefficient is positive, this follows from the Stroock--Varadhan support
theorem (e.g., \cite{IkedaWat}, Chapter VI, Theorem 8.1). Hence $%
P_{t_2}(x_\diamondsuit,(v^{\prime\prime}, v^{\prime\prime}))>0,$
and the
required statement holds true with
\[
c_{u,v,u^{\prime},v^{\prime},T}=c_1 \inf_{x\in
[u,v]}P_{t_2}\bigl(x,
\bigl(u^{\prime\prime}, v^{\prime\prime}\bigr)\bigr)>0.
\]
%

\subsection{\texorpdfstring{Proof of Proposition \protect\ref{p32}}{Proof of Proposition 3.2}}\label{s42}

\textit{Statement 1.} Take, analogously to (\ref
{phidef}), a
function $\psi:\mathbb{X}\to[1,+\infty)$ of the form
\[
\psi=\psi_\lozenge+\psi_\blacklozenge,
\]
where $\psi_\lozenge,\psi_\blacklozenge\in C^2(0,\infty)$, $\psi_\lozenge=0$
on $[2,\infty)$, $\psi_\blacklozenge=0$ on $(0,1]$,
\[
\psi_\lozenge(x)=x^{-\gamma^{\prime}}\qquad\mbox{for $x$ small enough,}\qquad
\psi_\blacklozenge(x)=x^{\delta^{\prime}}\qquad\mbox{for $x$ large enough,}
\]
with
\[
\gamma^{\prime}\in \biggl((\gamma-1)\vee0, {\frac{\alpha
}{2}}-1 \biggr),
\qquad \delta^{\prime}\in \biggl(\delta, {\frac{\beta}{2}} \biggr).
\]
Then, by the statement 1 of Proposition \ref{p21}, $\psi\in
\Dom(\mathcal{A})$
and $\psi$ satisfies (\ref{lyap}). By (\ref{xlarge}), one has
\[
\mathcal{A}\psi(x)\sim- C_\infty x^{\delta^{\prime}} = - C_\infty
\bigl(\phi(x)%
\bigr)^{\delta^{\prime}/\delta},\qquad x\to\infty
\]
with
\[
C_\infty= \theta\delta^{\prime} \biggl(1-{\frac{\delta^{\prime
}-1}{\beta
/2-1}}%
\biggr)>0.
\]
By (\ref{xsmall}), one has
\[
\mathcal{A}\psi(x)\sim- C_0 x^{-\gamma^{\prime}-1}= - C_0
\bigl(\phi (x) \bigr)%
^{(\gamma^{\prime}+1)/\gamma},\qquad x\to0
\]
with
\[
C_0=\theta\gamma^{\prime}\kappa \biggl(1-{
\frac{\gamma^{\prime
}+1}{\alpha
/2}}%
\biggr)>0.
\]
Finally, for every segment $[u,v]\subset(0,\infty)$ and every
$\varepsilon>0
$ one has
\[
\sup_{x\in[u,v]}{\phi(x)}<\infty,\qquad \sup_{x\in[u,v]}{
\frac
{|\mathcal{A}
\psi(x)|}{\phi^{1+\varepsilon}(x)}}<\infty,
\]
because $\phi, \mathcal{A} \psi\in C(0,\infty)$ and $\phi\geq1$. These
observations provide (\ref{est}) with small enough $c^{\prime},
\varepsilon$
and large enough $C^{\prime}$.

\textit{Statement 2.} By the elementary inequality $(\sum_{k=1}^m
a_k)^{1+\eps}\leq m^\eps\sum_{k=1}^m a_k^{1+\eps}$, we have
%
\begin{equation}
\label{53} \int_\XX \bigl(\phi(x_1)+\cdots+
\phi(x_m) \bigr)^{1+\varepsilon}\mu^t_{t_1,
\ldots, t_m}(
\mathrm{d}x)\leq m^\eps\sum_{k=1}^m
\int_\XX\phi^{1+\eps} \,\mathrm{d}
\mu^t_{t_k}.
\end{equation}

By the definition of $\mathcal{A}$, we have for arbitrary $\mu\in\Pf$
\[
E_{\mu}\psi(X_{t})=E_{\mu}\psi(X_{0})+E_{\mu}
\int_{0}^{t}\mathcal {A}%
\psi(X_{s})\,\mathrm{d}s.
\]
Together with (\ref{est}), this yields
%
\begin{eqnarray}
\label{54} %
\int_\XX\phi^{1+\eps} \,
\mathrm{d}\mu^t&=&\frac{1}{ t}\int_0^tE_\mu
\phi^{1+\eps
}(X_s)\,\mathrm{d}s\leq\frac{1}{ c't}E_\mu
\biggl[\int_0^tC' \,
\mathrm{d}s-\int_0^t\Af\psi
(X_s)\,\mathrm{d}s \biggr]
\nonumber
\\
&=&\frac{C'}{ c'}+\frac{1}{
c't}E_\mu\psi(X_0)-
\frac{1}{
c't}E_\mu\psi(X_t)
\\
&\leq&\frac{C'}{ c'}+\frac{1}{ c't}E_\mu\psi
(X_0)=\frac{C'}{ c'}+\frac{1}{ c't}\int
_\XX\psi \,\mathrm{d}\mu;
\nonumber
\end{eqnarray}
in the second inequality, we have used that $\psi$ is nonnegative. By
Corollary \ref{c31} with $\psi$ instead of $\phi$, we have
\[
\int_\XX\psi \,\mathrm{d}\mu_{t_k}\leq
\frac{C}{ c}+\int_\XX\psi \,\mathrm{d}\mu\leq \biggl(
\frac{C}{ c}+1 \biggr)\int_\XX\psi \,\mathrm{d}\mu,
\qquad k=1, \ldots,m
\]
because $\psi\geq1$. Using (\ref{53}) and (\ref{54}) with $\mu_{t_k},
k=1, \ldots,m$ instead of $\mu$, we obtain (\ref{Cesarobounds}).

\subsection{\texorpdfstring{Proof of Theorem \protect\ref{t31}}{Proof of Theorem 3.1}}\label{s43}

\textit{Statement 1.}
In \cite{Kulik}, Theorem 2.1, it is proved that a Markov process $X$
admits an exponential $\phi$-coupling under the following assumptions:
\begin{enumerate}[(iii)]
\item [(i)] $\phi\in \Dom(\Af)$ and (\ref{lyap}) holds true;

\item [(ii)] every \emph{level set} $\{\phi\leq R\}, R\geq1$ has a compact
closure in $\XX$;

\item [(iii)] for every compact $K\subset\XX$ there exists $T>0$ such that
%
\begin{equation}
\label{d} \sup_{x,x'\in K}\bigl\|P_T(x,\cdot)-P_T
\bigl(x',\cdot\bigr)\bigr\|_{\mathrm{var}}<2,
\end{equation}
where $\|\cdot\|_{\mathrm{var}}$ denotes the total variation norm.
\end{enumerate}

In our setting, (i) and (iii) are provided by Proposition \ref{p21}
(statements 1 and 2, resp.). Assumption (ii) holds true
trivially because $\phi(x)\to+\infty$ when either $x\to0$ or $x\to
\infty$. Hence, the required statement follows by Theorem 2.1 in \cite{Kulik}.
%
\begin{rem} In \cite{Kulik}, the notion of an exponential $\phi
$-coupling was introduced in a form, slightly weaker than the one from
Definition \ref{d13}; see the discussion after Definition \ref{d13}.
One can see easily that the proof of Theorem 2.1 in \cite{Kulik} can be
extended straightforwardly to provide an exponential $\phi$-coupling in
the sense of Definition \ref{d13}.
\end{rem}

\textit{Statement 2.} By statement 1, for a given $\mu\in\Pf$ there
exists a $(\mu,\pi)$-coupling which satisfies (\ref{phi-coup}). From
this fact, we will deduce (\ref{41}). In a particular case $\phi
\equiv
1, m=1$ such an implication is well known, and the proof for general
$\phi, m$ does not require any substantial changes when compared with
the standard one. To keep the exposition self-sufficient, we explain
the argument briefly.
Denote $\varkappa_t=\mu_{t+t_1, \ldots, t+t_m}-\pi_{t_1, \ldots, t_m}$,
\begin{eqnarray*}
\nu_{i,t}(\mathrm{d}y)&=&P \bigl(\bigl(Z_{t_1+t}^i,
\ldots,Z_{t_m+t}^i\bigr)\in \mathrm{d}y,
\\
\bigl(Z_{t_1+t}^1,\ldots,Z_{t_m+t}^1
\bigr)&\not=&\bigl(Z_{t_1+t}^2,\ldots ,Z_{t_m+t}^2
\bigr) \bigr),\qquad i=1,2.
\end{eqnarray*}
For arbitrary measurable function $f:\XX^m\to[0,+\infty)$, one has
%
\begin{eqnarray}
\label{42} %
\int_{\XX^m} f \,\mathrm{d}
\varkappa_t&=&Ef\bigl(Z_{t_1+t}^1,\ldots
,Z_{t_m+t}^1\bigr)-Ef\bigl(Z_{t_1+t}^2,
\ldots,Z_{t_m+t}^2\bigr)
\nonumber
\\[-8pt]
\\[-8pt]
&=&\int_{\XX^m} f \,\mathrm{d}\nu_{1,t}-\int
_{\XX^m} f \,\mathrm{d}\nu_{2,t}\leq\int
_{\XX^m} f \,\mathrm{d}\nu_{1,t}.
\nonumber
\end{eqnarray}
Denote by $A_{t}^+$ a set such that $\varkappa_t^+$ is supported by
$A_{t}^+$ and $\varkappa_t^-(A_{t}^+)=0$. By (\ref{42}), we have for
any measurable $A\subset A_{t}^+$:
\[
\varkappa_t^+(A)=\varkappa_t(A)\leq\nu_{1,t}(A).
\]
Because $\varkappa_t^+$ is supported by $A_{t}^+$, this gives finally
\[
\varkappa_t^+\leq\nu_{1,t}.
\]
Similarly,
\[
\varkappa_t^-\leq\nu_{2,t}.
\]
From these inequalities, we have
\begin{eqnarray*}
\|\varkappa_t\|_{\phi,\mathrm{var}}&\leq&\int_{\XX^m}
\bigl(\phi(x_1)+\cdots +\phi (x_m) \bigr)
\nu_{1,t}(\mathrm{d}x)+\int_{\XX^m} \bigl(
\phi(x_1)+\cdots+\phi (x_m) \bigr) \nu_{1,t}(
\mathrm{d}x)
\\
&=&E \Biggl(\sum_{j=1}^m \bigl[\phi
\bigl(Z^1_{t+t_j}\bigr)+\phi \bigl(Z^2_{t+t_j}
\bigr) \bigr] \Biggr)\1_{(Z_{t_1+t}^1,\ldots,Z_{t_m+t}^1)\not
=(Z_{t_1+t}^1,\ldots,Z_{t_m+t}^1)}
\\
&\leq&\sum_{j=1}^m E \bigl[\phi
\bigl(Z^1_{t+t_j}\bigr)+\phi\bigl(Z^2_{t+t_j}
\bigr) \bigr]\1_{Z^1_{t+t_j}\not=Z^2_{t+t_j}} \leq m C \mathrm{e}^{-c t}\int
_\XX\phi \,\mathrm{d}\mu ,
\end{eqnarray*}
where the last inequality comes from the assumption (\ref{phi-coup}).

\textit{Statement 3.} Estimate (\ref{41}) with $m=1$ provides similar and
weaker estimate with $\|\cdot\|_{\mathrm{var}}$ instead of $\|\cdot\|_{\phi
,\mathrm{var}}$. It is another standard observation that such an estimate,
together with an estimate of the form
%
\begin{equation}
\label{phi-bound} \int_\XX\phi \,\mathrm{d}
\mu_t\leq\tilde C\int_\XX\phi \,\mathrm{d}\mu,
\qquad \mu\in\Pf , t\geq0,
\end{equation}
provide (\ref{mixing}). Again, we explain this argument briefly.

The $\sigma$-algebra $\Ff^X_{\geq r}$ is generated by the algebra
$\Ff^{X, cyl}_{\geq r}$ of the sets of the form
%
\begin{equation}
\label{B-cyl} B=\bigl\{\bigl(X(v_1), \ldots, X(v_m)\bigr)
\in C\bigr\},\qquad v_1, \ldots, v_m\geq r, C\in\Bf\bigl(
\XX^m\bigr), m\geq1.
\end{equation}
Hence, in the identity (\ref{beta-m}), we can replace $\sup_{B\in\Ff
^X_{\geq t+s}}$ by $\sup_{B\in\Ff^{X, cyl}_{\geq t+s}}$. On the other
hand, for every $B$ of the form (\ref{B-cyl}) with $r=t+s$, we have
\[
P_\mu\bigl(B|\Ff_s^X\bigr)=T_tf(X_s),
\qquad P_\mu(B)=\int_\XX T_{t+s}f \,
\mathrm{d}\mu
\]
with
\[
f(x)=P_x\bigl(\bigl(X(v_1-t-s), \ldots,
X(v_m-t-s)\bigr)\in C\bigr), \qquad x\in\XX
\]
and
\[
T_rf(x)=\int_\XX f(y)P_r(x,
\mathrm{d}y)=E_xf(X_r),
\]
the usual notation for the semigroup generated by the Markov process
$X$. We have
\begin{eqnarray*}
\bigl|P_\mu(B|\Ff_s)-P_\mu(B)\bigr|&\leq&\biggl
\llvert T_tf(X_s)-\int_\XX f \,
\mathrm{d}\pi \biggr\rrvert + \biggl\llvert \int_\XX f \,
\mathrm{d}\pi-\int_\XX T_{t+s}f \,\mathrm{d}\mu
\biggr\rrvert
\\
&\leq&\bigl\| P_t(X_s,\cdot)-\pi\bigr\|_{\mathrm{var}}+\|
\mu_{t+s}-\pi\|_{\mathrm{var}},
\end{eqnarray*}
here we have used that $\|f\|\leq1$. Therefore, we have
%
\begin{equation}
\label{43} \beta^\mu(t)\leq\sup_{s\geq0} \bigl(\|
\mu_{t+s}-\pi\|_{\mathrm{var}}+E_\mu \bigl\|
P_t(X_s,\cdot)-\pi\bigr\|_{\mathrm{var}} \bigr).
\end{equation}
Note that (the weaker version of) (\ref{41}) gives
\[
\|\mu_{t+s}-\pi\|_{\mathrm{var}}\leq C\mathrm{e}^{-ct}\int
_\XX\phi \,\mathrm{d}\mu,\qquad \bigl\| P_t(X_s,
\cdot)-\pi\bigr\|_{\mathrm{var}}\leq C\mathrm{e}^{-ct}\phi(X_s).
\]
These observations combined with (\ref{phi-bound}) provide (\ref
{mixing}) with $C'=C(1+\tilde C)$.

Recall that $\phi$ satisfies a condition of the form (\ref{lyap});
denote respective constants by $c_L,C_L$. Then Corollary \ref{c31}
yields (\ref{phi-bound}) with $\tilde C=\frac{C_L}{ c_L}+1$ because it is
supposed that $\phi\geq1$. These observations finally lead to (\ref{mixing})
with
%
\begin{equation}
\label{Cpr} C'=C \biggl(2+\frac{C_L}{ c_L} \biggr).
\end{equation}

\subsection{\texorpdfstring{Proof of Theorem \protect\ref{t37}}{Proof of Theorem 3.2}}\label{s431}

Let $\gamma',\delta'$ be the values introduced in the construction of
the function $\psi$, see Section \ref{s42}. Denote
\[
\lambda= \biggl(\frac{\gamma'}{\gamma}\wedge\frac{\delta'}{\delta
} \biggr)^{-1}.
\]
For any signed measure $\varkappa$ on $\Bf(\XX^m)$, by the H\"older
inequality, we have
\[
\|\varkappa\|_{\phi,\mathrm{var}}\leq \Biggl(\int_{\XX^m} \Biggl(
\sum_{j=1}^m\phi (x_j)
\Biggr)^{\sigma p} |\varkappa|(\mathrm{d}x) \Biggr)^{1/p} \Biggl(\int
_{\XX^m} \Biggl(\sum_{j=1}^m
\phi(x_j) \Biggr)^{(1-\sigma
) q} |\varkappa|(\mathrm{d}x)
\Biggr)^{1/q}
\]
for any $\sigma>0$ and any $p,q>1$ with $1/p+1/q=1$. We put
$p=(\lambda
\sigma)^{-1}$ and take $\sigma$ close enough to 0, so that $p>1$. Then
$\phi^{\sigma p}=\phi^{1/\lambda}$, and
\[
\phi^{1/\lambda}(x)=x^{-\gamma ((\gamma'/\gamma)\wedge(\delta
'/\delta) )}\leq x^{-\gamma(\gamma'/\gamma)}=\psi(x)
\]
for $x$ small enough,
\[
\phi^{1/\lambda}(x)=x^{\delta ((\gamma'/\gamma)\wedge(\delta
'/\delta) )}\leq x^{\delta(\delta'/\delta)}=\psi(x)
\]
for $x$ large enough. Because $\phi$ is continuous and $\psi\geq1$,
this means that
%
\begin{equation}
\label{46} \Biggl(\sum_{j=1}^m
\phi(x_j) \Biggr)^{\sigma p}\leq C \sum
_{j=1}^m\psi(x_j)
\end{equation}
with some constant $C$. We have
\[
\frac{1}{ q}=1-\lambda\sigma,\qquad (1-\sigma) q=\frac{1-\sigma}{1-\lambda
\sigma},
\]
and in the above construction $\sigma$ can be taken close enough to 0
in order to provide inequality $(1-\sigma) q\leq1+\eps.$ Then we
obtain, finally,
%
\begin{equation}
\label{52} \|\varkappa\|_{\phi,\mathrm{var}}\leq C\| \varkappa \|_{\psi,\mathrm{var}}^{1/p}
\Biggl(\int_{\XX^m} \Biggl(\sum_{j=1}^m
\phi (x_j) \Biggr)^{1+\eps} |\varkappa|(\mathrm{d}x)
\Biggr)^{1/q}.
\end{equation}

Because the weighted total variation norm is a norm indeed, we have
\begin{eqnarray*}
&&\biggl\llVert \int_{0}^T (
\mu_{t+t_1, \ldots, t+t_m}-\pi_{t_1, \ldots, t_m}) \,\mathrm{d}t\biggr\rrVert_{\phi, \mathrm{var}}
\\
&&\quad\leq\sum_{k=0}^{[T]-1} \biggl\llVert
\int_{k}^{k+1} (\mu_{t+t_1, \ldots, t+t_m}-
\pi_{t_1, \ldots, t_m})\,\mathrm{d}t\biggr\rrVert_{\phi,
\mathrm{var}}
\\
&&\qquad{} + \biggl\llVert \int_{[T]}^T (
\mu_{t+t_1, \ldots, t+t_m}-\pi_{t_1, \ldots, t_m}) \,\mathrm{d}t\biggr\rrVert_{\phi, \mathrm{var}}
\\
&&\quad =\sum_{k=0}^{[T]-1} \bigl\llVert (
\mu_{k} )^1_{t_1,\ldots, t_m}-\pi_{t_1,\ldots, t_m}\bigr
\rrVert_{\phi, \mathrm{var}}+\bigl(T-[T]\bigr) \bigl\llVert (\mu_{[T]}
)^{T-[T]}_{t_1,\ldots, t_m}-\pi_{t_1,\ldots,
t_m}\bigr\rrVert_{\phi, \mathrm{var}};
\end{eqnarray*}
recall that $\mu_t$ denotes the one-dimensional distribution, see
(\ref
{mut}), and $\mu^t_{t_1, \ldots, t_m}$ denotes the Ces\`{a}ro mean, see
(\ref{cesaro}). By (\ref{52}), we have
\begin{eqnarray*}
&&\bigl\llVert (\mu_{k} )^1_{t_1,\ldots, t_m}-
\pi_{t_1,\ldots,
t_m}\bigr\rrVert_{\phi, \mathrm{var}}
\\
&&\quad\leq C \bigl\llVert (\mu_{k} )^1_{t_1,\ldots,
t_m}-
\pi_{t_1,\ldots, t_m}\bigr\rrVert_{\psi, \mathrm{var}}^{1/p}\bigl\llVert (
\mu_{k} )^1_{t_1,\ldots, t_m}-\pi_{t_1,\ldots, t_m}\bigr
\rrVert_{\phi
^{1+\eps
}, \mathrm{var}}^{1/q}
\\
&&\quad\leq C \bigl\llVert (\mu_{k} )^1_{t_1,\ldots,
t_m}-
\pi_{t_1,\ldots, t_m}\bigr\rrVert_{\psi, \mathrm{var}}^{1/p}
\\
&&\qquad{}\times \Biggl(\int_{\XX^m} \Biggl(\sum
_{j=1}^m\phi(x_j)
\Biggr)^{1+\eps} \bigl[ (\mu_{k} )^1_{t_1,\ldots, t_m}+
\pi_{t_1,\ldots, t_m} \bigr](\mathrm{d}x) \Biggr)^{1/q}.
\end{eqnarray*}

Recall that $\psi$ satisfies conditions of Proposition \ref{p21}. In
addition, it has compact level sets; see condition (ii) in Section \ref
{s43}. Then (\ref{41}) with $\psi$ instead of $\phi$ holds true, and
we have
\begin{eqnarray*}
\bigl\llVert (\mu_{k} )^1_{t_1,\ldots, t_m}-
\pi_{t_1,\ldots,
t_m}\bigr\rrVert_{\psi, \mathrm{var}}^{1/p} &=&\biggl\llVert
\int_{k}^{k+1}(\mu_{t_1+t,\ldots,
t_m+t}-
\pi_{t_1,\ldots, t_m})\,\mathrm{d}t\biggr\rrVert_{\psi, \mathrm{var}}^{1/p}
\\
& \leq& \biggl(\int_{k}^{k+1}\|
\mu_{t_1+t,\ldots, t_m+t}-\pi_{t_1,\ldots,
t_m}\|_{\psi, \mathrm{var}}\, \mathrm{d}t
\biggr)^{1/p}
\\
&\leq& m^{1/p}C^{1/p}\mathrm{e}^{-ck/p} \biggl(\int
_\XX \psi \,\mathrm{d}\mu \biggr)^{1/p}
\end{eqnarray*}
with the constants $c, C$ from (\ref{41}). Note that $\phi^{1+\eps}$ is
integrable w.r.t. $\pi$; see Remark \ref{r31}. Then
\begin{eqnarray*}
\int_{\XX^m} \Biggl(\sum_{j=1}^m
\phi(x_j) \Biggr)^{1+\eps} \pi_{t_1,\ldots
, t_m}(\mathrm{d}x)&
\leq& m^\eps\int_{\XX^m} \sum
_{j=1}^m\phi^{1+\eps}(x_j)
\pi_{t_1,\ldots, t_m}(\mathrm{d}x)
\\
&=&m^{1+\eps}\int_{\XX}\phi^{1+\eps} \,
\mathrm{d}\pi <\infty .
\end{eqnarray*}
On the other hand, by (\ref{Cesarobounds}) with $t=1$ we have
\begin{eqnarray*}
\int_{\XX^m} \Biggl(\sum_{j=1}^m
\phi(x_j) \Biggr)^{1+\eps} (\mu_{k}
)^1_{t_1,\ldots, t_m}(\mathrm{d}x)&=&\int_{\XX^m}
\Biggl(\sum_{j=1}^m\phi
(x_j) \Biggr)^{1+\eps} \mu^1_{t_1+k,\ldots, t_m+k}(
\mathrm{d}x)
\\
&\leq& C\int_\XX \psi \,\mathrm{d}\mu.
\end{eqnarray*}
Using the elementary inequality
\[
(x+y)^{1/q}\leq x^{1/q}+y^{1/q},\qquad x,y>0, q>1
\]
and the assumption $\psi\geq1$, we get from the above estimates
%
\begin{equation}
\label{47} \bigl\llVert (\mu_{k} )^1_{t_1,\ldots, t_m}-
\pi_{t_1,\ldots,
t_m}\bigr\rrVert_{\phi, \mathrm{var}}\leq\tilde C_m
\mathrm{e}^{-ck/p}\int_\XX\psi \,\mathrm{d}\mu
\end{equation}
with some explicitly calculable $\tilde C_m$. Similarly to (\ref{47})
(we omit the details), one can show that
%
\begin{equation}
\label{471}\bigl(T-[T]\bigr) \bigl\llVert (\mu_{[T]}
)^{T-[T]}_{t_1,\ldots, t_m}-\pi_{t_1,\ldots,
t_m}\bigr\rrVert_{\phi, \mathrm{var}}
\leq\tilde C_m \mathrm{e}^{-c[T]/p}\int_\XX
\psi \,\mathrm{d}\mu.
\end{equation}
From (\ref{47}) and (\ref{471}), we obtain the required inequality with
$C_m=\tilde C_m \sum_{k=0}^\infty \mathrm{e}^{-ck/p}$.

\subsection{\texorpdfstring{Proof of Theorem \protect\ref{t32}}{Proof of Theorem 3.3}}\label{s44}

In order to simplify
the notation, we assume $k=1$ and remove respective subscripts, that
is, write $f,\gamma, \delta$ instead of $f_i, \gamma_i, \delta_i$. One
can see that the proof below can be extended to the multidimensional
case easily; to do that, it is enough to replace the one-dimensional
``deviation inequalities'' (\ref{di1}) and (\ref{di2}) by completely
analogous inequalities for the components $f_i, i=1, \ldots, k$ of the
multidimensional function $f$.

We proceed in two steps: the ``coupling'' one and the ``truncation'' one.

The ``coupling'' step deals with the case where for some positive $\eps
$ the initial distribution $\mu$ satisfies (\ref{mu}).
Let $\phi$ be defined by (\ref{phidef}) with $\gamma, \delta$ from
(\ref{fi}). Then Theorem \ref{t31} provides that there exists a
$(\mu
,\pi)$-coupling $(Z^1,Z^2)$ for the process $X$, which satisfies (\ref
{phi-coup}). We have
\begin{eqnarray*}
E_\mu\Biggl\llvert \frac{1}{ n}\sum
_{l=1}^nf (X_{t_1+l}, \ldots,
X_{t_r+l} )- a_f\Biggr\rrvert &=&E\Biggl\llvert
\frac{1}{ n}\sum_{l=1}^nf
\bigl(Z^1_{t_1+l}, \ldots, Z^1_{t_r+l}
\bigr)- a_f\Biggr\rrvert
\\
&\leq& E\Biggl\llvert \frac{1}{ n}\sum_{l=1}^nf
\bigl(Z^2_{t_1+l}, \ldots, Z^2_{t_r+l}
\bigr)- a_f\Biggr\rrvert
\\
&&{}+\frac{1}{ n}\sum_{l=1}^n E
\bigl\llvert f \bigl(Z^1_{t_1+l}, \ldots, Z^1_{t_r+l}
\bigr)- f \bigl(Z^2_{t_1+l}, \ldots, Z^2_{t_r+l}
\bigr)\bigr\rrvert
\end{eqnarray*}
because $Z^2$ has the same distribution with $\{X^{st}(t), t\geq0\}$.
Recall that $X$ is ergodic, see \cite{Genon-Catalot}. Then, by the
Birkhoff--Khinchin theorem,
\[
E\Biggl\llvert \frac{1}{ n}\sum_{l=1}^nf
\bigl(X_{t_1+l}^{st}, \ldots, X_{t_r+l}^{st}
\bigr)- a_f\Biggr\rrvert \to0,\qquad n\to\infty.
\]
On the other hand, by (\ref{fi}) we have
\begin{eqnarray*}
&&E\bigl\llvert f \bigl(Z^1_{t_1+l}, \ldots,
Z^1_{t_r+l} \bigr)- f \bigl(Z^2_{t_1+l},
\ldots, Z^2_{t_r+l} \bigr)\bigr\rrvert
\\
&&\quad\leq C \sum_{j=1}^r E \bigl(\phi
\bigl(Z^1_{t_j+l}\bigr)+\phi\bigl(Z^2_{t_j+l}
\bigr) \bigr)\1_{(Z^1_{t_1+l}, \ldots,
Z^1_{t_r+l})\not=(Z^2_{t_1+l}, \ldots, Z^2_{t_r+l})}
\\
&&\quad \leq C \sum_{j=1}^r\sum
_{i=1}^r E \bigl(\phi\bigl(Z^1_{t_j+l}
\bigr)+\phi \bigl(Z^2_{t_j+l}\bigr) \bigr)\1_{Z^1_{t_i+l}\not=Z^2_{t_i+l}}
\end{eqnarray*}
(note that $C$ here does not coincide with the constant $C$ in (\ref
{fi}) because $\phi(x)\not=x^{-\gamma}+x^\delta$). By the H\"older
inequality and the elementary inequality $(a+b)^p\leq
2^{p-1}(a^p+b^p)$, $a,b\geq0, p>1$, we have for arbitrary $p,q>1$ with
$1/p+1/q=1$
\begin{eqnarray*}
&&E \bigl(\phi\bigl(Z^1_{t_j+l}\bigr)+\phi
\bigl(Z^2_{t_j+l}\bigr) \bigr)\1_{Z^1_{t_i+l}\not
=Z^2_{t_i+l}}
\\
&&\quad \leq2^{(p-1)/p} \bigl(E \bigl(\phi^p\bigl(Z^1_{t_j+l}
\bigr)+\phi^p\bigl(Z^2_{t_j+l}\bigr) \bigr)
\bigr)^{1/p} \bigl(P\bigl({Z^1_{t_i+l}\not
=Z^2_{t_i+l}}\bigr) \bigr)^{1/q}.
\end{eqnarray*}
We can take $p>1$ close enough to 1, so that $\gamma'=\gamma p<\gamma
+\eps, \delta'=\delta p<\delta+\eps$, and $\gamma', \delta'$ satisfy
(\ref{gammadelta}). Then $\phi'=\phi^p$ clearly has the form (\ref
{phidef}) with $\gamma', \delta'$ instead of $\gamma, \delta$.
Corollary \ref{c31} applied to $\phi'$ instead of $\phi$ yields that
\[
\sup_{t\geq0}E\phi^p\bigl(Z^1_{t}
\bigr)<\infty,\qquad \sup_{t\geq0}E\phi^p\bigl(Z^2_{t}
\bigr)=\int_\XX\phi^p \,\mathrm{d}\pi<\infty.
\]
On the other hand, (\ref{phi-coup}) and standing assumption $\phi\geq
1$ yield
\[
P\bigl(Z^1_{t}\not=Z^2_{t}\bigr)
\leq C \mathrm{e}^{-c t}\int_\XX\phi \,\mathrm{d}
\mu, \qquad t\geq0,
\]
where $c, C$ are the same as in (\ref{phi-coup}). Summarizing all the
above, we obtain
%
\begin{equation}
\label{di1} E\bigl\llvert f \bigl(Z^1_{t_1+l}, \ldots,
Z^1_{t_r+l} \bigr)- f \bigl(Z^2_{t_1+l},
\ldots, Z^2_{t_r+l} \bigr)\bigr\rrvert \leq C'
\sum_{i=1}^r \mathrm{e}^{-c(t_i+l)/q}
\end{equation}
with the same constant $c$ and some constant $C'$ which depends on
$\phi
, p, \mu$, and the constants $C$ in (\ref{fi}) and (\ref
{phi-coup}). Therefore
\[
E_\mu\Biggl\llvert \frac{1}{ n}\sum
_{l=1}^nf (X_{t_1+l}, \ldots,
X_{t_r+l} )- a_f\Biggr\rrvert \to0,\qquad n\to\infty,
\]
which completes the proof of statement 1 under the assumption (\ref
{mu}). To prove statement 2, we need to show that for any bounded
Lipschitz continuous function $F\dvtx \Re\to\Re$
%
\begin{equation}
\label{44} E_\mu F\bigl(S_n(X)\bigr)\to\int
_\Re F(y)\nu_f(\mathrm{d}y),
\end{equation}
where $\nu_f\sim\mathcal{N}(0, \Sigma_f^d)$ and
\[
S_n(X)=\frac{1}{\sqrt n}\sum_{l=1}^n
\bigl(f (X_{t_1+l}, \ldots, X_{t_r+l} )-a_f \bigr).
\]

In \cite{ALS}, Remark 3.1, it was shown that the general result by
Genon-Catalot \emph{et al.} (see \cite{Genon-Catalot}, Corollary 2.1) can be
applied to prove that the stationary Fisher--Snedecor diffusion is an
$\alpha$-mixing process with an exponential decay rate. Then the CLT
for $\alpha$-mixing sequences (see \cite{Hall}) provide
\[
E F\bigl(S_n\bigl(X^{st}\bigr)\bigr)\to\int
_\Re F(y)\nu_f(\mathrm{d}y).
\]
On the other hand, the estimates similar to those made above provide that
%
\begin{eqnarray}
\label{di2} %
&&\bigl|E_\mu F\bigl(S_n(X)\bigr)-E F
\bigl(S_n\bigl(X^{st}\bigr)\bigr)\bigr|
\nonumber\\[-8pt]\\[-8pt]
&&\quad\leq\frac{\Lip(F)}{\sqrt
n} \sum_{l=1}^nE \bigl|f
\bigl(Z_{t_1+l}^1, \ldots, Z_{t_r+l}^1
\bigr)-f \bigl(Z_{t_1+l}^2, \ldots, Z_{t_r+l}^2
\bigr) \bigr|\leq\frac{C'\Lip(F)}{
\sqrt n}\quad
\nonumber
\end{eqnarray}
with some constant $C'$. This proves statement 2 under the assumption
(\ref{mu}).

The ``truncation'' step removes the assumption (\ref{mu}). For an
arbitrary $\mu$ and any $a\in(0,1)$ there exist $\mu_a, \mu^a\in
\Pf$
such that $\mu_a$ is supported in some segment $[u,v]\subset(0,
\infty
)$, and
\[
\mu=(1-a)\mu_a+a\mu^a.
\]
Then $P_\mu=(1-a) P_{\mu_a}+a P_{\mu^a}, $ and $\mu_a$ satisfies
(\ref
{mu}). Hence, for any $\zeta>0$
\begin{eqnarray*}
&&\operatorname{\lim\sup}\limits_{n\to\infty}P_{\mu} \Biggl(\Biggl\llvert
\frac{1}{ n}\sum_{l=1}^nf
(X_{t_1+l}, \ldots, X_{t_r+l} )- a_f\Biggr\rrvert >
\zeta \Biggr)
\\
&&\quad\leq a \operatorname{\lim\sup}\limits_{n\to\infty}P_{\mu^a} \Biggl(\Biggl\llvert
\frac{1}{ n}\sum_{l=1}^nf
(X_{t_1+l}, \ldots, X_{t_r+l} )- a_f\Biggr\rrvert >
\zeta \Biggr)\leq a.
\end{eqnarray*}
Because $a$ is arbitrary, this proves statement 1 for arbitrary $\mu$.
Similar argument proves (\ref{44}) for arbitrary $\mu$, and completes
the proof of the theorem.

\subsection{\texorpdfstring{Proof of Theorem \protect\ref{t34}}{Proof of Theorem 3.4}}\label{s46}

Again, we assume
$k=1$. We note that both statement 1 and statement 2 hold true under
the respective conditions of Theorem \ref{t32}. The proof of this fact
is analogous to the proof of Theorem~\ref{t32} and therefore is
omitted. The only difference is that, in this proof, one requires the
continuous-time version of the CLT (\ref{CLTc}) for the stationary
version $X^{st}$ of the process $X$ instead of the discrete-time one.
This statement can be easily derived from the respective discrete-time
one by the standard discretization argument (see, e.g., \cite
{Billingsley}, pages 178--179). Hence, our task is to reduce the
conditions of Theorem \ref{t32} to those of Theorem \ref{t34}.

First, note that we can increase slightly $\gamma$, so that the
conditions of Theorem \ref{t34} still hold true. Let $\phi$ be defined
by (\ref{phidef}) with this new $\gamma$ and $\delta$ from the
formulation of the theorem. Because $\alpha>2$, condition (\ref
{gammadelta2}) yields (\ref{gammadelta1}). Then we can apply
Proposition \ref{p32} and define respective function $\psi$, see
Section \ref{s42}. While doing that, we can choose $\gamma', \delta'$
larger than, but close enough to $(\gamma-1)\vee0,\delta$,
respectively, so that
$\int_\XX\psi \,\mathrm{d}\mu<\infty$ if $\mu$ is supposed to satisfy (\ref{mu1}) and
%
\begin{equation}
\label{gammadelta5} \gamma'+\gamma<\frac{\alpha}{2},\qquad
\delta'+\delta<\frac{\beta}{2}
\end{equation}
if $\gamma, \delta$ satisfy (\ref{gammadelta2}). We put
\[
\|f\|_{\phi}=\sup_{x=(x_1, \ldots, x_r)}\frac{|f(x)|}{\sum_{j=1}^r\phi
(x_j)},\qquad
f_n(x)=f(x)\prod_{j=1}^r
\1_{x_j\geq1/n},\qquad n\geq1.
\]
For arbitrary $t_1, \ldots , t_r\geq0$ one has
\[
E \sum_{j=1}^r\phi\bigl(X^{st}_{t_j}
\bigr)=r\int_0^\infty\phi(x)\pi (\mathrm{d}x)<
\infty
\]
because $\gamma, \delta$ satisfy (\ref{gammadelta1}). Then, by (\ref
{fi}) and the Lebesgue dominated convergence theorem, $a_{f_n}\to a_f$.

We put $\tilde f_n=f_n+a_{f}-a_{f_n}$. Then the condition (\ref{fi})
with the initial $\gamma$ provide that
%
\begin{equation}
\label{fapr}\|f-\tilde f_n\|_{\phi}\to0,\qquad n\to\infty.
\end{equation}
On the other hand, $a_{\tilde f_n}=a_f$, and every $\tilde f_n$ satisfy
conditions of Theorem \ref{t32}. Hence, for every $n$
\begin{eqnarray*}
&&\operatorname{\lim\sup}\limits_{T\to\infty}E_\mu\biggl\llvert \frac{1}{ T}\int
_0^Tf (X_{t_1+t}, \ldots,
X_{t_r+t} )\,\mathrm{d}t- a_f\biggr\rrvert
\\
&&\quad\leq\operatorname{\lim\sup}\limits_{T\to\infty
}\frac{1}{ T}E_\mu\int
_0^T \bigl|f (X_{t_1+t}, \ldots,
X_{t_r+t} )-\tilde f_n (X_{t_1+t}, \ldots,
X_{t_r+t} ) \bigr| \,\mathrm{d}t
\\
&&\quad\leq \operatorname{\lim\sup}\limits_{T\to\infty}\frac{C\|f-\tilde f_n\|_{\phi}}{ T}E_\mu\int
_0^T\sum_{j=1}^r
\phi(X_{t_j+t})\,\mathrm{d}t
\\
&&\quad= C\|f-\tilde f_n\|_{\phi} \operatorname{\lim\sup}\limits_{T\to
\infty}\sum
_{j=1}^r \biggl(\frac{T+t_j}{ T}\int
_\XX\phi \,\mathrm{d}\mu^{T+t_j}-\frac{t_j}{ T}
\int_\XX\phi \,\mathrm{d}\mu^{t_j} \biggr).
\end{eqnarray*}
Then from (\ref{Cesarobounds}) with $m=1$ and $\eps=0$ we obtain that,
when $\mu$ satisfies (\ref{mu1}),
\[
\operatorname{\lim\sup}\limits_{T\to\infty}E_\mu\biggl\llvert \frac{1}{ T}\int
_0^Tf (X_{t_1+t}, \ldots,
X_{t_r+t} )\,\mathrm{d}t- a_f\biggr\rrvert \leq C \|f-
\tilde f_n\|_{\phi}
\]
with some constant $C$. Because $n$ is arbitrary and (\ref{fapr})
holds, this proves (\ref{LLNc}) in the mean sense. If (\ref{mu1})
fails, then (\ref{LLNc}) still holds in the sense of convergence of
probability; one can show this using the truncation argument from the
previous section. This proves statement 1.

Denote $Q=\max_j t_j-\min_jt_j$ and assume that $T>Q$. Then
\begin{eqnarray*}
&&E_\mu \biggl[\frac{1}{\sqrt T}\int_0^T
\bigl(f (X_{t_1+t}, \ldots, X_{t_r+t} )-\tilde f_n
(X_{t_1+t}, \ldots, X_{t_r+t} ) \bigr)\, \mathrm{d}t \biggr]^2
\\
&&\quad\leq\frac{2}{ T} \biggl[\int_0^T\int
_{s}^{T\wedge
(s+Q)}+\int_0^{T-Q}
\int_{s+Q}^T \biggr]E_\mu \bigl(f
(X_{t_1+t}, \ldots, X_{t_r+t} )-\tilde f_n
(X_{t_1+t}, \ldots, X_{t_r+t} ) \bigr)
\\
&&\hspace*{159pt}{}\times \bigl(f (X_{t_1+s}, \ldots, X_{t_r+s} )-\tilde
f_n (X_{t_1+s}, \ldots, X_{t_r+s} ) \bigr)\,
\mathrm{d}t\, \mathrm{d}s\\
&&\quad=:I_1+I_2.
\end{eqnarray*}
We estimate $I_1, I_2$ separately. We explain the estimates in the
particular case $r=2, t_1=0, t_2=Q$; the general case is quite
analogous, but the calculations are more cumbersome. We have
%
\begin{equation}
\label{49} I_1\leq\frac{C\|f-\tilde f_n\|_{\phi}^2}{ T}\int_0^T
\int_{s}^{T\wedge
(s+Q)} E_\mu \bigl(
\phi(X_{t})+\phi(X_{t+Q}) \bigr) \bigl(\phi
(X_{s})+\phi (X_{s+Q}) \bigr)\,\mathrm{d}t\, \mathrm{d}s.\quad
\end{equation}
By the Markov property of the process $X$,
\[
\int_0^T\int_{s}^{T\wedge(s+Q)}
E_\mu\phi(X_{t})\phi(X_{s})\,\mathrm{d}t \,\mathrm{d}s
\leq E_\mu\int_0^T
\phi(X_{s}) \biggl(\int_0^Q
T_v\phi(X_s) \,\mathrm{d}v \biggr)\,\mathrm{d}s;
\]
here we have used the standard notation
\[
T_vf(x)=\int_\XX f(y)P_t(x,
\mathrm{d}y).
\]
Note that $P_t(x,\cdot)= (\delta_x )_t$. Hence, by (\ref
{Cesarobounds}) with $m=1, \eps=0,$ and $\mu=\delta_x$, we have
%
\begin{equation}
\label{55} \int_0^Q T_v\phi(x)
\,\mathrm{d}v\leq QC\psi(x),\qquad x\in\XX.
\end{equation}
By the inequalities (\ref{gammadelta5}), the function $\Phi=\phi\psi$
has the form (\ref{phidef}) with the parameters satisfying (\ref
{gammadelta1}). Then, using once again (\ref{Cesarobounds}) with
$\Phi
$ instead of $\phi$, we get
\[
\int_0^T\int_{s}^{T\wedge(s+Q)}
E_\mu\phi(X_{t})\phi(X_{s})\,\mathrm{d}t \,\mathrm{d}s
\leq QC E_\mu\int_0^T
\Phi(X_{s})\,\mathrm{d}s\leq TQC'\int
_\XX\Phi \,\mathrm{d}\mu;
\]
the constants $C, C'$ here depend on $\phi, \psi$, etc., but does not
depend on $Q,T,$ and $\mu$. Similar calculations provide estimates for
other parts of the integral in the right-hand side of (\ref{49}). For
instance, changing the variables $s'=s+Q$ and using the Markov property
at the point $t\leq s'$, we get
\begin{eqnarray*}
&&\int_0^T\int_{s}^{T\wedge(s+Q)}
E_\mu\phi(X_{t})\phi(X_{s+Q})\,\mathrm{d}t \,\mathrm{d}s \\
&&\quad=
E_\mu\int_0^T
\phi(X_{t}) \biggl(\int_{Q\vee t}^{t+Q}
\phi(X_{s'}) \,\mathrm{d}s' \biggr)\,\mathrm{d}t
\\
&&\quad\leq E_\mu\int_0^T
\phi(X_{t}) \biggl(\int_0^Q
T_v\phi (X_t)\,\mathrm{d}v \biggr)\,\mathrm{d}t \leq
TQC'\int_\XX\Phi \,\mathrm{d}\mu;
\end{eqnarray*}
in the last inequality we use (\ref{55}) and (\ref{Cesarobounds}) with
$\Phi$ instead of $\phi$.

Summarising these estimates, we get
\[
I_1\leq CQ\|f-\tilde f_n\|_{\phi}^2
\int_\XX\Phi \,\mathrm{d}\mu.
\]

To estimate $I_2$, we use the Markov property at the time moment $s+Q$
and write
\[
I_2\leq\frac{C\|f-\tilde f_n\|_{\phi}}{ T}E_\mu\int_0^{T}
\bigl(\phi (X_s)+\phi(X_{s+Q}) \bigr)F^{n,Q,T}_{s}(X_{s+Q})
\,\mathrm{d}s
\]
with
\[
F^{n,Q,T}_{s}(x)=\biggl\llvert \int_{0}^{T-s-Q}
E_x \bigl(f (X_{t}, X_{t+Q} )-\tilde
f_n (X_{t}, X_{t+Q} ) \bigr)\,\mathrm{d}t\biggr
\rrvert .
\]
Denote $g_n=f-\tilde f_n$. Because, by the construction, $a_f=a_{\tilde
f_n}$, we have $\int_{\XX^2} g_n\,\mathrm{d}\pi_{t, t+Q}=0$ for every $t$. Then
\[
F^{n,Q,T}_{s}(x)= \biggl\llvert \int_\XX
g_n \,\mathrm{d} \biggl(\int_{0}^{T-s-Q} \bigl( (
\delta_x)_{t, t+Q} -\pi_{t, t+Q} \bigr)\,\mathrm{d}t
\biggr)\biggr\rrvert .
\]
Clearly,
\[
\biggl\llvert \int_{\XX^m} g \,\mathrm{d}\varkappa\biggr\rrvert \leq\|g
\|_{\phi}\| \varkappa\|_{\phi, \mathrm{var}}
\]
for any measurable function $g$ on $\XX^m$ and any signed measure
$\varkappa$. Then, by (\ref{45}),
\begingroup
\abovedisplayskip=6.7pt
\belowdisplayskip=6.7pt
\[
F^{n,Q,T}_{s}(x)\leq C\|f-\tilde f_n
\|_{\phi}\psi(x).
\]
Recall that $\psi$ satisfies the Lyapunov-type condition (\ref{lyap}).
Then by the Markov property and the moment bound from Corollary \ref
{c31} we have $E_\mu\phi(X_s)\psi(X_{s+Q})\leq C E_\mu\phi
(X_s)\psi
(X_s)$, which together with the preceding estimate gives
\[
I_2\leq\frac{C\|f-\tilde f_n\|_{\phi}^2}{ T}E_\mu\int_0^{T}
\bigl(\phi (X_s)\psi(X_s)+\phi(X_{s+Q})
\psi(X_{s+Q}) \bigr)\,\mathrm{d}s.
\]
Using once again (\ref{Cesarobounds}) with $\Phi=\phi\psi$ instead of
$\phi$ and recalling the estimates for $I_1$, we get finally
%
\begin{equation}
\label{48} E_\mu \biggl[\frac{1}{\sqrt T}\int_0^T
\bigl(f(X_t)-\tilde f_n(X_t) \bigr) \,\mathrm{d}t
\biggr]^2\leq C \|f-\tilde f_n\|^2_\phi
\int_\XX\Phi \,\mathrm{d}\mu.
\end{equation}

By the construction, every $f_n$ satisfies conditions of Theorem \ref
{t32}, and therefore (\ref{CLTc}) holds true with $f_n$ instead of
$f$. Then, if $\Phi$ is integrable w.r.t. $\mu$, (\ref{48}) and the
approximation argument, similar to the one used in the proof of Theorem
\ref{t32}, lead to (\ref{CLTc}) for $f$ with
%
\begin{equation}
\label{sigmaf1} \Sigma_f^c=\lim_{n\to\infty}
\Sigma_{f_n}^c.
\end{equation}
On the other hand, if we write
\[
\Sigma_{f, R}^c=\int_{-R}^R
\Cov \bigl(f \bigl(X_{t_1+t}^{st}, \ldots, X_{t_r+t}^{st}
\bigr), f \bigl(X_{t_1}^{st}, \ldots, X_{t_r}^{st}
\bigr)\bigr) \,\mathrm{d}t,
\]
then
%
\begin{equation}
\label{481} \bigl|\Sigma_{f, R}^c-\Sigma_{\tilde f_n, R}^c\bigr|
\leq C \|f-\tilde f_n\|^2_\phi\int
_\XX\Phi \,\mathrm{d}\pi;
\end{equation}
the proof of (\ref{481}) is similar to the proof of (\ref{48}) and is
omitted. Therefore the integral (\ref{sigmaf}) coincides with the limit
(\ref{sigmaf1}). This completes the proof of statement 2 when $\Phi$ is
integrable w.r.t. $\mu$. For general $\mu$, we use the truncation
argument from the previous section.

\subsection{\texorpdfstring{Proof of Theorem \protect\ref{tA}}{Proof of Theorem 3.5}}

Again, we restrict ourselves by
the case $k=1$. The proof is based on the following auxiliary estimate.

\begin{lem}\label{lA1} Under conditions of Theorem \ref{tA}, for any $T$
\[
E \biggl(\int_0^T \bigl(f
(X_{t_1+t}, \ldots, X_{t_r+t} )-a_f \bigr) \,\mathrm{d}t
\biggr)^2\leq CT \|f\|_\phi^2 \int
_\XX\Phi \,\mathrm{d}\mu
\]
with some $\Phi$ satisfying conditions of statement 1 of Proposition
\ref{p21}.
\end{lem}\eject
\endgroup
\begin{pf} We assume that $f$ is centered and $r=1$. The general
case can be reduced to this one using the same arguments with those
explained Section \ref{s46}.

We proceed like in Section \ref{s46}: take $\psi$ of the form (\ref
{phidef}) with $\gamma'\in((\gamma-1)\vee0, \alpha/2-1), \delta'<\beta/2$ such that $\gamma+\gamma'<\alpha/2-1, \delta+\delta'<\beta
/2$ and put $\Phi=\phi\psi$.
Then
\begin{eqnarray*}
E \biggl(\int_0^T f(X_t)\,
\mathrm{d}t \biggr)^2&=&2\int_0^T
Ef(X_s)\int_s^T
f(X_t)\,\mathrm{d}t \,\mathrm{d}s
\\
&\leq&2\|f\|_\phi\int_0^T
E\bigl|f(X_s)\bigr|\biggl\llVert \int_0^{T-s}
\bigl((\delta_{X_s})_{r}-\pi \bigr)\,\mathrm{d}r\biggr
\rrVert_{\phi,\mathrm{var}}\,\mathrm{d}s
\\
&\leq& C\|f\|_\phi\int_0^T
E\bigl|f(X_s)\bigr|\psi(X_s)\,\mathrm{d}s,
\end{eqnarray*}
here we have used the Markov property and Theorem \ref{t37}. On the
other hand, Corollary \ref{c31} applied to $\Phi$ instead of $\phi$ gives
\[
\int_0^T E\bigl|f(X_s)\bigr|
\psi(X_s)\,\mathrm{d}s \leq\|f\|_\phi\int
_0^T \biggl(\int_\XX\Phi
\,\mathrm{d}\mu_s \biggr)\,\mathrm{d}s\leq C\|f\|_\phi T
\int_\XX\Phi \,\mathrm{d}\mu %
\]
with some other constant $C$, which completes the proof.
\end{pf}

Let us proceed with the proof of the theorem. By Theorem \ref{t34},
finite-dimensional distributions of $Y_T$ converge to that of $B$.
Hence, we need to prove the weak compactness, only. In addition, it is
sufficient to prove weak compactness in $D([0,1])$ instead of
$C([0,1])$: when we succeed to do that, we get the weak convergence
$Y_T\Rightarrow B$ in $D([0,1])$. Because both $Y_T$ and $B$ have
continuous trajectories, this would imply the weak convergence
$Y_T\Rightarrow B$ in $C([0,1])$.

For the function $\Phi$ constructed in the proof of Lemma \ref{lA1},
there exists $q>1$ such that $\Phi^q$ still satisfies conditions of
Proposition \ref{p21}, statement 1. Then, for $p$ such that
$1/p+1/q=1$, we have for every $v_1<v_2<v_3$
%
\begin{eqnarray}
\label{A1} %
&&E\bigl|Y_T(v_1)-Y_T(v_2)\bigr|^{2/ p}\bigl|Y_T(v_2)-Y_T(v_3)\bigr|^2\nonumber
\\
&&\quad\leq C\|f\|_\phi^2(v_3-v_2)
E\bigl|Y_T(v_1)-Y_T(v_2)\bigr|^{2/ p}
\Phi\bigl(X(v_2T)\bigr)
\\
&&\quad\leq C\|f\|_\phi^2(v_3-v_2)
\bigl(E\bigl|Y_T(v_1)-Y_T(v_2)\bigr|^{2}
\bigr)^{1/p} \bigl(E\Phi^q\bigl(X(v_2T)\bigr)
\bigr)^{1/q}\nonumber
\\
&&\quad\leq C\|f\|_\phi^{2+2/p}(v_3-v_2)
(v_2-v_1)^{1/p} E\Phi\bigl(X(v_1T)
\bigr)^{1/p} E\Phi^q\bigl(X(v_2T)
\bigr)^{1/q}\nonumber
\\
&&\quad \leq C\|f\|_\phi^{2+2/p}(v_3-v_1)^{1+1/p}
\biggl(\int_\XX\Phi \,\mathrm{d}\mu \biggr)^{1/p}
\biggl(\int_\XX\Phi^q \,\mathrm{d}\mu
\biggr)^{1/q}.\nonumber %
\end{eqnarray}
Here we have used subsequently Lemma \ref{lA1}, the H\"older
inequality, Lemma \ref{lA1} again, and Corollary \ref{c31} with $\Phi,
\Phi^q$ instead of $\phi$. Theorem 15.6 in \cite{Billingsley} and
(\ref
{A1}) provide weak compactness in $D([0,1])$ of the family $\{X_T\}$.

\subsection{\texorpdfstring{Proof of Theorem \protect\ref{t36}}{Proof of Theorem 4.1}}

By Example \ref{e35}, under
the assumptions of Theorem \ref{t36}, for any fixed $t>0$ either
$(\overline{m}_{{-1},c}, \overline{m}_{{1},c},\allowbreak \overline{m}_{2,c},
\overline{R}_{c}(t))$ or $(\overline{m}_{{-1},d}, \overline{m}_{{1},d},
\overline{m}_{2,d}, \overline{R}_{d}(t))$ is an asymptotically normal
estimator of $({m}_{{-1},c},\allowbreak {m}_{{1},c}, {m}_{2,c}, R(t))$. Note that
the assumption $\alpha>2, \beta>8$ (in the continuous-time case) is
equivalent to
\[
\{-1, 1, 2\}\in \biggl(-\frac{\alpha}{4}-\frac{1}{2},
\frac{\beta}{4} \biggr),
\]
while the assumption $\alpha>4, \beta>8$ (in the discrete-time case) is
equivalent to
\[
\{-1, 1, 2\}\in \biggl(- \biggl(\frac{\alpha}{2}-1 \biggr)\wedge \biggl(
\frac{\alpha
}{4} \biggr), \frac{\beta}{4} \biggr).
\]

The invariant distribution density for the process $X$ can be written
in the form
%
\begin{equation}
\label{fsdensD} \mathfrak{p}(x) =\frac{1}{ x B(\alpha/2, \beta/2)} \biggl(\frac{\alpha
x}{\alpha x+\varrho}
\biggr)^{\alpha/2} \biggl(\frac{\varrho}{\alpha x+
\varrho} \biggr)^{\beta/2}
\end{equation}
with $\varrho={(\beta-2)\kappa/\beta}$. Respective moments are equal
%
\begin{equation}\label{momentseqD} %
 m_\upsilon=\int
_0^\infty x^\upsilon\mathfrak{p}(x)
\,\mathrm{d}x= \biggl(\frac{\varrho
}{\alpha} \biggr)^\upsilon\frac{\Gamma(\alpha/2+\upsilon) \Gamma
(\beta
/2-\upsilon)}{\Gamma(\alpha/2)\Gamma(\beta/2)},\qquad
\upsilon\in \biggl(-\frac{\alpha}{2}, \frac{\beta}{2} \biggr).
\end{equation}
In particular,
\[
m_{-1}=\frac{\alpha}{(\alpha-2)(\beta-2)\kappa},\qquad m_1=\frac{\kappa}{
\beta},\qquad
m_2=\frac{(\alpha+2)(\beta-2)\kappa^2}{\alpha(\beta
-4)\beta^2}.
\]
On the other hand, one has
\[
\Corr \bigl(X_0^{st}, X_t^{st}
\bigr)=\mathrm{e}^{-\theta t},
\]
see \cite{Bibby1}, Theorem 2.3(iii). Resolving the above identities
for a fixed $t$, we can write $(\alpha, \beta, \kappa, \theta
)=G(m_{-1}, m_1, m_2, R(t))$ with
\begin{eqnarray*}
G_1(x,y,z,w)&=&\frac{2( xyz-y^2)}{ xyz-2z+y^2 },\qquad G_2(x,y,w)=
\frac{4x(z-y^2)}{ xz-2xy^2+ y},
\\
G_3(x,y,z,w)&=&\frac{4xy(z-y^2)}{ xz-2xy^2+ y},\qquad G_4(x,y,z,w)=-
\frac{1}{
t}\log \biggl(\frac{w}{ z-y^2} \biggr).
\end{eqnarray*}
Clearly, the function $G$ is well defined and smooth in some
neighbourhood of the point
\[
\mathbf{x}= \bigl(m_{-1}(\alpha, \beta, \kappa, \theta),
m_1(\alpha , \beta , \kappa, \theta), m_2(\alpha, \beta,
\kappa, \theta), \bigl[R(t) \bigr](\alpha, \beta, \kappa, \theta) \bigr).
\]
Then one can obtain the required statements using the continuity
mapping theorem and the functional delta method (see \cite{Serfling},
Theorem 3.3.A). Asymptotic covariance matrices for $(\widehat{\alpha
}_c, \widehat{\beta}_c, \widehat{\kappa}_c, \widehat{\theta}_c)$ and
$(\widehat{\alpha}_d, \widehat{\beta}_d, \widehat{\kappa}_d,
\widehat
{\theta}_d),$ are given by the formula
%
\begin{equation}
\label{covar} \Sigma_c(\alpha, \beta,\kappa, \theta)=D
\Sigma_c D^{\top},\qquad \Sigma_d(\alpha, \beta,
\kappa, \theta)=D\Sigma_d D^{\top},
\end{equation}
where $\Sigma_c, \Sigma_d$ are the asymptotic covariance matrices for
\[
\bigl(\overline{m}_{{-1},c}, \overline{m}_{{1},c},
\overline{m}_{2,c}, \overline{R}_{c}(t)\bigr),\qquad \bigl(
\overline{m}_{{-1},d}, \overline {m}_{{1},d}, \overline{m}_{2,d},
\overline{R}_{d}(t)\bigr),
\]
respectively, and $D_{ij} =  [ \frac{\partial
G_{i}}{\partial x_{j}}  ](\mathbf{x})$, $i,j \in\{1, 2,3,4\}$.

\section*{Acknowledgements} The authors express their deep gratitude to
the referee for numerous comments and suggestions that were helpful and
led for substantial improvement of the paper.

A.M. Kulik was supported in part by the State fund for fundamental
researches of Ukraine and the Russian foundation for basic research,
Grant F40.1/023.

N.N. Leonenko was supported in part by grant of the European commission
PIRSES-GA-2008-230804 (Marie Curie).



%

\printhistory

\end{document}